\documentclass[12pt,reqno]{amsart}
\usepackage{amssymb,amsthm,amsmath,amstext,amsxtra}
\usepackage{mathrsfs,bm}
\usepackage{fullpage}
\usepackage[all]{xy}
\usepackage{booktabs}
\usepackage{mathtools}
\usepackage{hyperref}
\hypersetup{colorlinks=true,urlcolor=blue,citecolor=blue,linkcolor=blue}
\usepackage{enumerate}
\usepackage{ stmaryrd }
\usepackage{enumitem}
\usepackage{comment}
\usepackage{multirow}
\usepackage{colonequals}
\usepackage{tikz}
\usepackage{marginnote}
\usepackage[export]{adjustbox}

\usepackage{tikz}
\usepackage{tikz-cd}
\usepackage{xcolor}

\numberwithin{equation}{section}

\newtheorem{example}{\text{Example}}[section]
\newtheorem{theorem}{Theorem}[section]

\newtheorem{remark}[theorem]{Remark}
\newtheorem{definition}[theorem]{Definition}

\begin{document}

\title{K3 mirror symmetry, Legendre family and Deligne's conjecture for Fermat quartic}

\author{Wenzhe Yang}
\address{SITP, Physics Department, Stanford University, CA, 94305}
\email{yangwz@stanford.edu}

\begin{abstract}
In this paper, we will study the connections between the mirror symmetry of K3 surfaces and the geometry of the Legendre family of elliptic curves. We will prove that the mirror map of the Dwork family is equal to the period map of the Legendre family. This result provides an interesting explanation to the modularities of counting functions for K3 surfaces from the mirror symmetry point of view. We will also discuss the relations between the arithmetic geometry of smooth fibers of the Fermat pencil (Dwork family) and that of the smooth fibers of the Legendre family, e.g. Shioda-Inose structures, zeta functions, etc. In particular, we will study the relations between the Fermat quartic, which is modular with a weight-3 modular form $\eta(4z)^6$, and the elliptic curve over $\lambda=2$ of the Legendre family, whose weight-2 newform is labeled as \textbf{32.2.a.a} in LMFDB. We will also compute the Deligne's periods of the Fermat quartic, which are given by special values of the theta function $\theta_3$. Then we will numerically verify that they satisfy the predictions of Deligne's conjecture on the special values of $L$-functions of critical motives.
\end{abstract}

%\vspace{-30pt}
\maketitle
\setcounter{tocdepth}{1}
\vspace{-13pt}
\tableofcontents
\vspace{-13pt}

\section{Introduction}
The mirror symmetry of K3 surfaces is significantly different from that of Calabi-Yau threefolds, and it can be described in terms of pure Hodge structures on the lattice of total integral cohomology groups \cite{AM}. Given a K3 surface $X$, its total integral cohomology group 
\begin{equation}
H^*(X,\mathbb{Z})=H^0(X,\mathbb{Z}) \oplus H^2(X,\mathbb{Z}) \oplus H^4(X,\mathbb{Z})
\end{equation}
is a free $\mathbb{Z}$-module of rank 24 with a Mukai pairing that is even unimodular. Together with this Mukai pairing, $H^*(X,\mathbb{Z})$ is isomorphic to the enlarged K3 lattice $\Lambda \oplus U$, where $\Lambda$ is the K3 lattice
\begin{equation}
\Lambda=E_8(-1)^2 \oplus U^3.
\end{equation}
On the lattice $H^*(X,\mathbb{Z})$, there is a weight-two pure Hodge structure defined by \cite{Hartmann}
\begin{equation}
\begin{aligned}
H^{2,0}_B(X)&=H^{2,0}(X),\\
H^{1,1}_B(X)&=H^0(X,\mathbb{C}) \oplus H^{1,1}(X) \oplus H^4(X,\mathbb{C}),\\
H^{0,2}_B(X)&=H^{0,2}(X).\\
\end{aligned}
\end{equation}
Let us call the data $H_B(X,\mathbb{Z}):=\left(H^*(X,\mathbb{Z}), \langle \cdot \rangle,H^{p,q}_B(X) \right)$ the B-model pure Hodge structure of $X$ \cite{AM,Hartmann}. 

The A-model pure Hodge structure of a K3 surface depends on the choice of a K\"ahler form. Suppose $Y$ is a K3 surface with complex structure $I$ and K\"ahler form $\omega_Y$, then $\omega_Y$ defines a symplectic structure on the underlying differential manifold of $Y$. The A-model pure Hodge structure is associated to this symplectic manifold, which also depends on the choice of a $B$-field $\beta \in H^2(Y,\mathbb{R})$. The cohomology class $\mho$ is by definition
\begin{equation}
\mho=\exp \left( \beta+i \omega_Y \right)=1+\left( \beta+i \omega_Y \right)+\frac{1}{2}\left( \beta+i \omega_Y \right)^2 \in H^*(X,\mathbb{C}),
\end{equation}
and with respect to the Mukai pairing, it satisfies
\begin{equation}
\langle \mho,\mho \rangle=0,~ \langle \mho,\overline{\mho} \rangle >0.
\end{equation}
By demanding $H_A^{2,0}(Y)=\mathbb{C} \mho$, we get a pure Hodge structure on $H^*(Y,\mathbb{Z})$, and together with the Mukai pairing, we obtain the A-model pure Hodge structure, $H_A(Y,\mathbb{Z})$ \cite{Hartmann}. The K3 surface $X$ with a holomorphic twoform $\Omega$ and the K3 surface $Y$ with a complexified K\"ahler form $\mho$ are said to form a mirror pair if there exists a Hodge isometry \cite{AM, Hartmann}
\begin{equation} \label{eq:mirrorpairdefn}
H_A(Y,\mathbb{Z}) \simeq H_B(X,\mathbb{Z}).
\end{equation}
It should be noticed that this definition can be viewed as a further refinement of Dolgachev's work \cite{Dolgachev}.

In the paper \cite{Hartmann}, Hartmann has proved that smooth quartic K3 surfaces and the Dwork family form a mirror pair in the sense of \ref{eq:mirrorpairdefn}. More precisely, given a smooth quartic surface $Y \subset \mathbb{P}^3$, the Fubini-Study K\"ahler form on $\mathbb{P}^3$ induces a symplectic structure $\omega_Y$  on $Y$. By a result of Moser \cite{Moser}, all smooth quartics are symplectomorphic to each other. Now introduce a B-field $\tau_1 \omega_Y,\tau_1 \in \mathbb{R}$ and define a complexified K\"ahler form 
\begin{equation}
\tau \omega_Y=(\tau_1+i\tau_2)\omega_Y, \tau_1 \in \mathbb{R},\tau_2 \in \mathbb{R}^+.
\end{equation}
So on the K\"ahler side, we have obtained a family of complexfified symplectic manifolds over the upper half plane $\mathbb{H}$ of $\mathbb{C}$
\begin{equation} \label{eq:kahlerside}
\rho:\mathscr{Y} \rightarrow \mathbb{H},
\end{equation}
where the fiber at $\tau \in \mathbb{H}$ is $(Y,\mho=\exp \left(\tau \omega_Y \right))$. 

The mirror family of $\mathscr{Y}$ \ref{eq:kahlerside} can be constructed from the Fermat pencil of K3 surfaces
\begin{equation} \label{eq:introFermatpencil}
\mathscr{F}_\psi:\{ x_0^4+x_1^4+x_2^4+x_3^4-4\psi \,x_0x_1x_2x_3 =0 \} \subset \mathbb{P}^3
\end{equation} 
by first taking the quotient with respect to a $(\mathbb{Z}/4\mathbb{Z})^4$ action and then a minimal resolution of singularities. The result is a family of K3 surfaces over $\mathbb{P}^1$ with $\psi$ as its parameter
\begin{equation} \label{eq:introdworkfamily}
\pi:\mathscr{X} \rightarrow \mathbb{P}^1,
\end{equation} 
which is called the Dwork family. By abuse of notations, the underlying differential manifold of a K3 surface is also denoted by $X$. There exists a canonical method to construct a nowhere vanishing holomorphic twoform $\Omega^f_\psi$ on a smooth fiber $\mathscr{F}_\psi$ of \ref{eq:introFermatpencil}, which further induces a holomorphic twoform $\Omega_\psi$ on the smooth fiber $\mathscr{X}_\psi$ of the Dwork family \ref{eq:introdworkfamily}. From the papers \cite{Hartmann, Nagura}, there exist two integral homology cycles $\Gamma_1$ and $\Gamma_2 $ in $H_2(X,\mathbb{Z})$ such that
\begin{equation}
\begin{aligned}
W_0(\psi)&=\int_{\Gamma_1} \Omega_\psi=\sum_{n=0}^\infty \frac{(4n)!}{(n!)^4(4 \psi)^{4n}},\\
W_1(\psi)&=\int_{\Gamma_2} \Omega_\psi=\frac{1}{2 \pi i} \left(-4W_0 \cdot \log(4 \psi)+4 \sum_{n=0}^\infty \frac{(4n)!}{(n!)^4(4\psi)^{4n}}[\Psi(4n+1)-\Psi(n+1)] \right),
\end{aligned}
\end{equation} 
where $\Psi$ is the polygamma function. The mirror map between the families $\mathscr{Y}$ and $\mathscr{X}$ is given by the quotient \cite{Hartmann, Nagura}
\begin{equation} \label{eq:intromirrormap}
\tau=\frac{W_1(\psi)}{W_0(\psi)}.
\end{equation}
It has been shown in the paper \cite{Hartmann} that $\mathscr{Y}_\tau$ and $\mathscr{X}_\psi$ form a mirror pair under the mirror map \ref{eq:intromirrormap}. While it is expected in \cite{Hartmann} that this mirror map may play an important role in the homological mirror symmetry of quartic K3 surfaces \cite{Seidel}.

The motivation of this paper is to explore the interesting connections between the previously stated mirror symmetry of K3 surfaces and the geometry of the Legendre family of elliptic curves. One important result of this paper is that under a transformation of the form 
\begin{equation} \label{eq:introquadratictransformation}
\frac{1}{\psi^4}=\frac{16(1-\lambda)\lambda^2}{(\lambda-2)^4},
\end{equation}
the mirror map \ref{eq:intromirrormap} is the same as the period map of the Legendre family. Recall that the Legendre family of elliptic curve is defined by the equation \cite{CarlsonMullerPeters}
\begin{equation}
\mathscr{E}_\lambda: y^2=x(x-1)(x-\lambda),
\end{equation}
which has a nowhere vanishing holomorphic oneform $\omega_\lambda=dx/2y$. The underlying differential manifold of a smooth fiber is the torus $T=S^1 \times S^1$. There is a choice of a basis $\{\gamma_0,\gamma_1 \}$ for $H_1(T,\mathbb{Z})$ such that the periods of $\omega_\lambda$ are \cite{Chandra, YangSW}
\begin{equation}
\int_{\gamma_0} \omega_\lambda=2 \pi \, \varpi_0(\lambda)=2 \pi \,_2F_1(\frac{1}{2},\frac{1}{2};1;\lambda),~\int_{\gamma_1} \omega_\lambda=2 \pi \,\varpi_1(\lambda).
\end{equation} 
The period map of the Legendre family is given by the quotient
\begin{equation} \label{eq:introperiodLegendre}
\tau=\frac{\varpi_1(\lambda)}{\varpi_0(\lambda)},
\end{equation}
the inversion of which is the famous modular lambda function \cite{Chandra, YangSW}
\begin{equation}
\lambda(\tau)=16q-128 q^2+704q^3-3072 q^4+11488q^5-38400q^6+ \cdots; q=\exp (\pi i \tau).
\end{equation}
Using the theory of hypergeometric functions, we will show that under the transformation \ref{eq:introquadratictransformation} we have
\begin{equation} 
\begin{aligned}
W_0(\psi)&=(1-\frac{\lambda}{2})\, \varpi^2_0(\lambda),\\
W_1(\psi)&=(1-\frac{\lambda}{2})\, \varpi_0(\lambda)\varpi_1(\lambda).
\end{aligned}
\end{equation}
Hence the mirror map \ref{eq:intromirrormap} of the Dwork family is essentially the same as the period map \ref{eq:introperiodLegendre} of the Legendre family. Moreover, this property provides an interesting explanation to the modularities of counting functions for K3 surfaces from the mirror symmetry point of view.

In fact, the holomorphic twoform $\Omega_\psi$ of the Dwork family satisfies a Picard-Fuchs equation $\mathcal{D}_3 \Omega_\psi=0$ with $\mathcal{D}_3$ a third order differential operator. Moreover, $\mathcal{D}_3$ is the symmetric square of a second order differential operator $\mathcal{D}_2$, which has two linearly independent solutions of the form
\begin{equation}\label{eq:introperiodsLegendreR}
\begin{aligned}
\pi_0^K(\lambda)&=(1-\frac{\lambda}{2})^{1/2} \, \varpi_0(\lambda),\\
\pi_1^K(\lambda)&=(1-\frac{\lambda}{2})^{1/2} \, \varpi_1(\lambda).
\end{aligned}
\end{equation}
Therefore $(\pi_0^K)^2$, $\pi_0^K \pi_1^K$ and $(\pi_1^K)^2$ are linearly independent solutions of $\mathcal{D}_3$. Based on this property, we will show that the pure Hodge structure on the transcendental lattice of $\mathscr{X}_\psi$ is isomorphic to that on the transcendental lattice of $\mathscr{E}_\lambda \times \mathscr{E}_\lambda$. The Picard number of a smooth fiber $\mathscr{X}_\psi$ satisfies $\geq 19$, hence from Morrison's work \cite{Morrison}, $\mathscr{X}_\psi$  admits a Shioda-Inose structure. An interesting question is to look at the connections between the Shioda-Inose structure of $\mathscr{X}_\psi$ and the geometry of $\mathscr{E}_\lambda \times \mathscr{E}_\lambda$.

On the other hand, the holomorphic twoform $\Omega^f_\psi$ of the Fermat pencil satisfies the same Picard-Fuchs equation as $\Omega_\psi$, i.e. $\mathcal{D}_3 \Omega^f_\psi=0$. Hence the previous results for the periods of $\Omega_\psi$ also work for $\Omega^f_\psi$. In particular, the pure Hodge structure on the transcendental lattice of $\mathscr{F}_\psi$ is isomorphic to that on the transcendental lattice of $\mathscr{E}_\lambda \times \mathscr{E}_\lambda$. The Picard number of a smooth fiber $\mathscr{F}_\psi$ satisfies $\geq 19$, hence it is also very interesting to study the connections between the Shioda-Inose structure of $\mathscr{F}_\psi$ and the geometry of $\mathscr{E}_\lambda \times \mathscr{E}_\lambda$. We will also discuss the connections between the zeta functions of $\mathscr{E}_\lambda$ and that of $\mathscr{F}_\psi$. One important example is the smooth fiber $\mathscr{F}_0$, called the Fermat's quartic, 
\begin{equation} \label{eq:introFermatQ}
 x_0^4+x_1^4+x_2^4+x_3^4=0,
\end{equation}
which is modular and associated to it is a weight-3 modular form $\eta^6(4z)$. Under the transformation \ref{eq:introquadratictransformation}, the point $\psi=0$ corresponds to $\lambda=2$. The minimal Weierstrass equation of the elliptic curve $\mathscr{E}_2$ at $\lambda=2$ of the Legendre family is
\begin{equation}
y^2=x^3-x,
\end{equation}
which is labeled as \textbf{32.a3} in LMFDB. The weight-2 newform associated to $\mathscr{E}_2$ is labeled as \textbf{32.2.a.a} in LMFDB. We will see the modular form $\eta^6(4z) $ can be considered as the symmetric square of \textbf{32.2.a.a}.

Furthermore, the formula \ref{eq:introperiodsLegendreR} allows us to explicitly compute the Deligne's periods of the Fermat quartic \ref{eq:introFermatQ}. More precisely, let $\mathbf{M}_0$ be the two dimensional pure motive that corresponds to the transcendental cycles of the Fermat quartic. From Deligne's  paper \cite{DeligneL}, the Tate twist $\mathbf{M}_0 \otimes \mathbb{Q}(n)$ is critical if and only if $n=1,2$. Using the method developed in the paper \cite{YangDeligne}, we find that Deligne's periods $c^+(\mathbf{M}_0 \otimes \mathbb{Q}(n))$, $n=1,2$, are given by special values of the theta function $\theta_3$
\begin{equation}
\begin{aligned}
c^+(\mathbf{M}_0 \otimes \mathbb{Q}(1))&=(2 \pi i)\,\theta^4_3(0,-i e^{-\pi/2}),\\
c^+(\mathbf{M}_0 \otimes \mathbb{Q}(2))&=i(2 \pi i)^2 \,\theta^4_3(0,-i e^{-\pi/2}).
\end{aligned}
\end{equation}
The $L$-function associated to $\mathbf{M}_0$ is just $L(\eta^6(4z),s)$, which has an integral representation 
\begin{equation} \label{eq:intromellintransformL}
L(\eta(4z)^6,s)=\frac{(2 \pi)^s}{\Gamma(s)} \int_0^\infty \eta(4iz)^6z^s \frac{dz}{z},
\end{equation}
thus its special values at $s=1,2$ can be numerically evaluated. Using Mathematica, we will numerically show that
\begin{equation}
\begin{aligned}
c^+(\mathbf{M}_0 \otimes \mathbb{Q}(1))&=16 \,L(\mathbf{M}_0 \otimes \mathbb{Q}(1),0),\\
c^+(\mathbf{M}_0 \otimes \mathbb{Q}(2))&=-64\, L(\mathbf{M}_0 \otimes \mathbb{Q}(2),0),
\end{aligned}
\end{equation}
which indeed satisfy the predictions of Deligne's conjecture \cite{DeligneL,YangDeligne}.

The outline of this paper is as follows. In Section \ref{sec:K3surfacesgeometry}, we will give an overview of the geometry of K3 surfaces, which includes a short review of the Shioda-Inose structures of algebraic K3 surfaces. In Section \ref{sec:mirrorfamilyconstruction}, we will introduce the Fermat pencil and the construction of the Dwork family. We will introduce the solutions of the Picard-Fuchs equation of the Dwork family, and the construction of the mirror map. In Section \ref{sec:legendrefamilyofellipticcurves}, we will first review some elementary properties of elliptic curves defined over $\mathbb{Q}$. Then we will introduce the Legendre family of elliptic curves, its periods and the modular lambda function. In Section \ref{sec:mirrormapLegendre}, we will use the quadratic transformations of hypergeometric functions to show the mirror map of the Dwork family is the same as the period map of the Legendre family. We will also discuss the connections between this property and the modularity of counting functions for K3 surfaces. Section \ref{sec:ShiodaInose} is a brief discussion of the potential relations between the Shioda-Inose structures of the Dwork family (or the Fermat pencil) and the Legendre family. Section \ref{sec:zetafunctionsconnections} studies the relations between the zeta functions of smooth fibers of the Fermat pencil and that of the Legendre family. In Section \ref{sec:DeligneconjectureFermatquartic}, we will look at the relations between the weight-3 newform of the Fermat quartic and the weight-2 newform of the elliptic curve at $\lambda=2$ of the Legendre family. We will also explicitly compute the Deligne's periods for the Fermat quartic, and numerically verify that they satisfy the predictions of Deligne's conjecture. Section \ref{sec:conclusions} contains a summary of the conclusions of this paper and some further open questions. Appendix \ref{sec:Weilconjectures} is a short review of the Weil conjectures.

\section{The geometry of K3 surfaces} \label{sec:K3surfacesgeometry}

In this section, we will review some elementary geometric properties of K3 surfaces, e.g. N\'eron-Severi group and transcendental lattice, etc. We will also give a brief overview of the Shioda-Inose structures of algebraic K3 surfaces with Picard numbers $\geq 19$. The readers who are familiar with these elementary materials can simply skip this section.

\subsection{An overview of K3 surfaces}

A K3 surface is by definition a 2-dimensional complex manifold $X$ with vanishing first sheaf cohomology and trivial canonical bundle \cite{Huybrechts}
\begin{equation}
H^1(X,\mathscr{O}_X)=0,~\Omega_X= \mathscr{O}_X.
\end{equation}
Here $\mathscr{O}_X$ is the sheaf of holomorphic functions on $X$ and $\Omega_X$ is the canonical bundle, i.e. the sheaf of holomophic twoforms on $X$. From its definition, a K3 surface is a two dimensional Calabi-Yau manifold. The triviality of $\Omega_X$ immediately implies there exists a nowhere vanishing holomorphic twoform $\Omega$ on $X$. In the definition, we have included non-algebraic K3 surfaces, and in fact most K3 surfaces are non-algebraic \cite{Huybrechts}. But every K3 surface is K\"ahler, and any two K3 surfaces are deformation equivalent to each other, in particular they are diffeomorphic to each other \cite{Huybrechts}. 

\begin{remark}\label{remarkDFK}
Since every two K3 surfaces are diffeomorphic to each other, in this paper the symbol $X$ will also mean the underlying differential manifold of K3 surfaces.
\end{remark}

The integral cohomology groups of a K3 surface is torsion free \cite{Peters}, and in fact we have:
\begin{enumerate}
\item $H^1(X,\mathbb{Z})=H^3(X,\mathbb{Z})=0$.

\item $H^2(X,\mathbb{Z})$ is a lattice of rank 22.

\end{enumerate}
From Hodge theory, there exist Hodge decompositions on the cohomology groups of a K\"ahler manifold. For example, the Hodge decomposition on $H^2(X,\mathbb{Z})$ is of the form
\begin{equation} \label{eq:hodgedecompostionH2}
H^2(X,\mathbb{Z}) \otimes_{\mathbb{Z}} \mathbb{C}=H^{2,0}(X) \oplus H^{1,1}(X) \oplus H^{0,2}(X),
\end{equation}
which defines a pure Hodge structure on $H^2(X,\mathbb{Z})$. The Hodge number $h^{i,j}$ is by definition $\text{dim}\,H^{i,j}(X)$, and the Hodge diamond of a K3 surface is of the form
\begin{center} 
\begin{tabular}{ c c c c c  }
   &  & 1 &  &    \\ 
   & 0&   & 0&    \\   
  1&  & 20 &  & 1.   \\  
   & 0&   & 0&    \\   
   &  & 1 &  &    \\
\end{tabular}
\end{center}

The Picard group of $X$, denoted by $\text{Pic}(X)$, is the abelian group of isomorphism classes of line bundles on $X$ \cite{Huybrechts}. The first Chern class defines a homomorphism from it to $H^2(X,\mathbb{Z})$
\begin{equation}
c_1: \text{Pic}(X) \rightarrow H^2(X,\mathbb{Z}).
\end{equation}
From the Lefschetz theorem on $(1,1)$-classes, the image of $c_1$ is $H^{1,1}(X) \cap H^2(X,\mathbb{Z})$, which is called the N\'eron-Severi group $\text{NS}(X)$ of $X$ \cite{Huybrechts}. The group $\text{NS}(X)$ is also characterized by the property \cite{Peters,Huybrechts,Yui}
\begin{equation} \label{eq:neronseveriOmega}
\gamma \in \text{NS}(X) \iff \int_X \gamma \smile \Omega=0,
\end{equation}
where $\smile$ means the cup product between cohomological cycles. The group $\text{NS}(X)$ is a sub-lattice of $H^2(X,\mathbb{Z})$, whose rank, denoted by $\rho(X)$, is called the Picard number of $X$. As the Hodge number $h^{1,1}(X) \leq 20$, we deduce $\rho(X) \leq 20$. 

There is a cup-product pairing on $H^2(X,\mathbb{Z})$
\begin{equation} \label{eq:k3middlecupproduct}
\langle \alpha,\beta \rangle = \int_X \alpha \smile \beta;~\alpha,\beta \in H^2(X,\mathbb{Z}),
\end{equation}
which is even unimodular with signature $(3, 19)$. Together with this cup-product pairing, the lattice $H^2(X,\mathbb{Z})$ is isomorphic to \cite{Huybrechts}
\begin{equation}
E^2_8(-1) \oplus (U_2)^3.
\end{equation}
The cup-product pairing \ref{eq:k3middlecupproduct} induces a non-degenerate symmetric bilinear form on $H^2(X,\mathbb{R})$, whose restriction to the subspace 
\begin{equation}
H^{1,1}_{\mathbb{R}}(X):=H^{1,1}(X) \cap H^2(X,\mathbb{R})
\end{equation}
is a non-degenerate symmetric bilinear form with signature $(1,19)$. The cup-product pairing \ref{eq:k3middlecupproduct} induces a pairing on $\text{NS}(X)$ with signature $(1,\rho(X)-1)$. The transcendental lattice $T(X)$ is by definition the orthogonal complement of $\text{NS}(X)$ with respect to the cup-product pairing \ref{eq:k3middlecupproduct}
\begin{equation}
T(X):=\text{NS}(X)^\perp \subset H^2(X,\mathbb{Z}).
\end{equation}
Similarly, the cup-product pairing \ref{eq:k3middlecupproduct} induces an even unimodular pairing on the lattice $T(X)$. Moreover, the pure Hodge structure on $H^2(X,\mathbb{Z})$ induces pure Hodge structures on the two sub-lattices $\text{NS}(X)$ and $T(X)$. 

Under the cup-product pairing \ref{eq:k3middlecupproduct}, the holomorphic twoform $\Omega$ satisfies
\begin{equation}
\langle\Omega,\Omega \rangle=0,~ \langle \Omega, \overline{\Omega} \rangle >0.
\end{equation}
In fact, the pure Hodge structure on $H^2(X,\mathbb{Z})$ is completely determined by the twoform $\Omega$ since $H^{2,0}(X)$ (resp. $H^{0,2}(X)$) is spanned by $\Omega$ (resp. $\overline{\Omega}$) and 
\begin{equation}
H^{1,1}(X)=\left( H^{2,0}(X) \oplus H^{0,2}(X)\right)^\perp.
\end{equation}

More generally, on the total integral cohomology group of $X$
\begin{equation}
H^*(X,\mathbb{Z})=H^0(X,\mathbb{Z}) \oplus H^2(X,\mathbb{Z}) \oplus H^4(X,\mathbb{Z}),
\end{equation}
there exists a Mukai pairing defined by 
\begin{equation} \label{eq:mukaipairing}
\langle (a_0,a_2,a_4),(b_0,b_2,b_4) \rangle= \int_X \left( a_2 \smile b_2-a_0 \smile b_4 -a_4 \smile b_0 \right).
\end{equation}
The Mukai pairing is important in the study of the mirror symmetry of K3 surfaces, and the readers are referred to the papers \cite{AM, Dolgachev, Hartmann} for more details.

\subsection{The Shioda-Inose structures of algebraic K3 surfaces}

Before we discuss the Shioda-Inose structures of algebraic K3 surfaces, let us first recall the definition of Kummer surfaces.
\begin{definition}
Suppose $A$ is an abelian surface with involution $\iota$, then the quotient variety $A/\iota$ has 16 $A_1$ singularities that correspond to the 2-division points of $A$. The minimal resolution of $A/\iota$ is a K3 surface $\text{Km}(A)$ that is called Kummer surface. 
\end{definition}

Let us first discuss the Shioda-Inose structures of singular K3 surfaces. Given a complex K3 surface $X$, it is called singular if its Picard number is 20. Here singular does not mean it is not smooth, but means such a K3 surface is exceptional. In many ways, singular K3 surfaces behave like elliptic curves with complex multiplication (CM). The transcendental lattice $T(X)$ of a singular K3 surface $X$ is of rank-2, and the cup-product pairing $Q(X)$ on $T(X)$ is an even integral positive definite binary quadratic form
\begin{equation} \label{eq:quadraticform}
Q(X)=
\left(
\begin{array}{cc}
 2 a & b \\
 b & 2 c \\
\end{array}
\right);~a,b,c \in \mathbb{Z}.
\end{equation}
The determinant $d$ of $Q(X)$ is by definition $d=b^2-4ac$. From the works \cite{PS,Shioda}, the map $X \mapsto Q(X)$ is a bijection between the isomorphic classes of singular K3 surfaces and the even integral positive definite binary quadratic forms up to conjugations by elements of $\text{SL}_2(\mathbb{Z})$.

Now we briefly review the construction of the inverse map of $X \mapsto Q(X)$. Given an even integral positive definite binary quadratic form $Q$ of the form \ref{eq:quadraticform}, there are two isogenous elliptic curves $\mathcal{E}_\tau$ and $\mathcal{E}_{\tau'}$ with 
\begin{equation}
\tau=\frac{-b+\sqrt{d}}{2a},~\tau'=\frac{b+\sqrt{d}}{2},
\end{equation}
both of which admit CM in the field $\mathbb{Q}(\sqrt{d})$. Here $\mathcal{E}_\tau$ means the complex torus $\mathbb{C}/(\mathbb{Z}+\tau \mathbb{Z})$, etc. However, it turns out that the Kummer surface of $\mathcal{E}_\tau \times \mathcal{E}_{\tau'}$ is a singular K3 surface with intersection form $2Q$. To cure this defect, Shioda-Inose construct a special elliptic fibration for $\text{Km}(\mathcal{E}_\tau \times \mathcal{E}_{\tau'})$. Then they show there exists a suitable quadratic base change of this fibration, the pull-back with respect to which is a singular K3 surface $X$ with intersection form $Q$.  From this construction, every singular K3 surface is defined over a number field \cite{PS,Shioda}. In conclusion, we have a diagram of the form
\begin{equation}
\begin{tikzcd}[column sep=scriptsize]
X \arrow[dr] & & \mathcal{E}_\tau \times \mathcal{E}_{\tau'} \arrow[dl] \\
& \text{Km}(\mathcal{E}_\tau \times \mathcal{E}_{\tau'}) & 
\end{tikzcd},
\end{equation}
where the arrows are of degree 2 \cite{Morrison}. The map from $X$ to $\text{Km}(\mathcal{E}_\tau \times \mathcal{E}_{\tau'})$ in this diagram is a Nikulin involution. Recall that an involution $\iota$ on a K3 surface $X$ is called a Nikulin involution if it preserves the twoform $\Omega$, i.e.
\begin{equation}
\iota^*(\Omega)=\Omega.
\end{equation}
From \cite{Nikulin}, every Nikulin involution has eight isolated fixed points. Now let us introduce the concept of Hodge isometry. 
\begin{definition}
Suppose $A$ and $A'$ are two lattices endowed with pure Hodge structures and bilinear forms, then a Hodge isometry $A \rightarrow A'$ is an isomorphism that respects both the pure Hodge structures and the bilinear forms.
\end{definition}

A general K3 surface $X$ is said to admit a Shioda-Inose structure if there exists a Nikulin involution on $X$ and a diagram of rational maps
\begin{equation}
\begin{tikzcd}[column sep=scriptsize]
X \arrow[dr,dotted] & & A \arrow[dl,dotted] \\
& \text{Km}(A) & 
\end{tikzcd},
\end{equation}
where $A$ is an abelian surface. Here the dotted arrow from $X$ to $ \text{Km}(A)$ corresponds to the quotient by Nikulin involution, and both dotted arrows are rational maps of degree 2. Moreover, this diagram induces a Hodge isometry between the transcendental lattices
\begin{equation}
T(X)\simeq T(A).
\end{equation}
From Theorem 6.3 of \cite{Morrison}, a K3 surface $X$ admits a Shioda-Inose structure if and only if there exists an abelian surface $A$ and a Hodge isometry $T(X) \simeq T(A)$. While from Corollary 6.4 of \cite{Morrison}, algebraic K3 surfaces with Picard numbers $\geq 19$ always admit Shioda-Inose structures. The readers could consult \cite{Morrison} for more details about the Shioda-Inose structures of algebraic K3 surfaces.

\section{The Dwork family and its mirror map} \label{sec:mirrorfamilyconstruction}

In this section, we will first discuss the Picard-Fuchs equation of the Fermat pencil of K3 surfaces and its independent solutions. Then we will briefly review the Greeene-Plesser construction of the mirror family of  quartic K3 surfaces from the Fermat pencil, which is usually called the Dwork family of K3 surfaces. We will also look at the construction of the mirror map of the Dwork family and its properties, which have been studied in \cite{Nagura}.

\subsection{The Fermat pencil of K3 surfaces} \label{sec:canonicalperiodsoffermatpencil}

The adjuction formula tells us that a smooth quartic surface in $\mathbb{P}^3$ is K3 \cite{Huybrechts}.  The Fermat pencil of K3 surfaces is a pencil of quartic surfaces in $\mathbb{P}^3$ defined by
\begin{equation} \label{eq:fermatpencilK3}
\mathscr{F}_\psi:\{f_\psi=0 \} \subset \mathbb{P}^3, ~f_\psi:=x_0^4+x_1^4+x_2^4+x_3^4-4\psi \,x_0x_1x_2x_3,
\end{equation}
where $(x_0,x_1,x_2,x_3)$ form the projective coordinate of $\mathbb{P}^3$. In a more formal language, formula \ref{eq:fermatpencilK3} defines a family
\begin{equation}\label{eq:descendfermatpencil}
\pi^f: \mathscr{F} \rightarrow \mathbb{P}^1,
\end{equation}
which is in fact defined over $\mathbb{Q}$. The fiber $\mathscr{F}_\psi$ is smooth if and only if $\psi$ does not lie in
\begin{equation}
\Sigma=\{\psi^4=1 \} \cup \{ \infty \}.
\end{equation}
When $\psi^4=1$, the fiber $\mathscr{F}_\psi$ has 16 singularities of type $A_1$. While when $\psi=\infty$, the Fermat pencil degenerates into the union of four complex planes
\begin{equation}
x_0 x_1x_2x_3=0.
\end{equation}
The Picard number $\rho(\mathscr{F}_\psi) $ of a smooth fiber $\mathscr{F}_\psi$ is $\geq 19$ \cite{Huybrechts}. There exists a projective linear transformation
\begin{equation} \label{eq:equivariantTF}
x_0 \rightarrow \zeta_4 \,x_0, ~x_i \rightarrow x_i,~i=1,2,3;~\zeta_4=\exp \pi i/2,
\end{equation}
that induces an isomorphism between $\mathscr{F}_\psi$ and $\mathscr{F}_{\zeta_4 \psi}$. So the `true' parameter for the Fermat pencil \ref{eq:fermatpencilK3} is in fact the variable $t$ defined by
\begin{equation} \label{eq:tpsirelation}
t:=1/\psi^4.
\end{equation}

\begin{remark}
If $\psi \neq 0,\infty$, the following isomorphism
\begin{equation} \label{eq:equivariantTFpsi}
x_0 \rightarrow \frac{1}{\psi} \,x_0, ~x_i \rightarrow x_i,~i=1,2,3,
\end{equation}
transforms the fiber $\mathscr{F}_\psi$ to the rationally defined surface
\begin{equation} \label{eq:rationalfpsi}
t x_0^4+x_1^4+x_2^4+x_3^4-4x_0x_1x_2x_3=0.
\end{equation}
\end{remark}

On the smooth fiber $\mathscr{F}_\psi$ of the Fermat pencil \ref{eq:fermatpencilK3}, i.e. $\psi \neq \Sigma$, there is a canonical way to construct a nowhere vanishing holomorphic twoform $\Omega^f_\psi$ \cite{MarkGross,Nagura}. On $\mathbb{P}^3$, there is a meromorphic threeform $\Theta_\psi$ given by
\begin{equation} \label{eq:meromorphicthreeform}
\Theta_\psi=\sum_{i=0}^3 (-1)^i \frac{\psi\, x_i \,dx_0 \wedge \cdots \wedge \widehat{dx_i} \wedge \cdots \wedge d x_3}{f_\psi},
\end{equation}
which is a well-defined threeform on $\mathbb{P}^3-\mathscr{F}_\psi$. $\Theta_\psi$ is automatically closed, hence its residue along $\mathscr{F}_\psi$ is well-defined, which is by definition the holomorphic twoform $\Omega^f_\psi$
\begin{equation} \label{eq:quartictwoform}
\Omega^f_\psi:= \text{Res}(\Theta_\psi).
\end{equation}
More explicitly, take the open subset of $\mathscr{F}_\psi$ defined by $x_3=1$, then the residue of $\Theta_\psi$ is equal to \cite{MarkGross}
\begin{equation} \label{eq:meroHoloTwo}
\Omega^f_\psi=\psi \frac{dx_0 \wedge dx_1 }{\partial f_\psi /\partial x_2} \Big|_{\mathscr{F}_\psi}.
\end{equation}
In fact, it can be explicitly shown that the meromorphic twoform \ref{eq:meroHoloTwo}, which is a priori only defined on $\partial f_\psi /\partial x_2 \neq 0$, extends to a global nowhere vanishing twoform on $\mathscr{F}_\psi$  \cite{MarkGross}. Notice that for a rational $\psi$, $\Omega^f_\psi$ is defined over $\mathbb{Q}$. Moreover, for $\psi \in \Sigma$, the previous construction defines a twoform $\Omega^f_\psi$ that is nowhere vanishing on the smooth locus of $\mathscr{F}_\psi$.

On the fiber $\mathscr{F}_\psi$, there exists a homology cycle $\beta_0 \in H_2(X,\mathbb{Z})$ consists of the points \cite{MarkGross}
\begin{equation} \label{eq:beta0cycle}
|x_0|=|x_1|=\delta,~x_3=1,
\end{equation}
and $x_2$ given by the solution to $f_\psi=0$ that tends to 0 when $\psi \rightarrow \infty$. Notice that $\beta_0$ is a torus in $\mathscr{F}_\psi$ that is a continuous deformation of \cite{MarkGross}
\begin{equation} 
\{(x_0,x_1,x_2,x_3) \in \mathbb{P}^3:|x_0|=|x_1|=\delta,~x_2=0,~x_3=1 \} \subset \mathscr{F}_\infty.
\end{equation}
For large $\psi$, the integration of $\Omega^f_\psi$ over $\beta_0$, up to a nonzero rational multiple, has a power series expansion of the form \cite{MarkGross,Nagura}
\begin{equation}
(2 \pi i)^2\sum_{n=0}^\infty \frac{(4n)!}{(n!)^4(4 \psi)^{4n}},
\end{equation}
which converges in a neighborhood of $\psi=\infty$. However it is practically very difficult to find other periods of $\Omega^f_\psi$ by explicit integration, that is where the Picard-Fuchs equation comes to the rescue.

\subsection{The Picard-Fuchs equation}

The Griffiths transversality tells us that the holomorphic twoform $\Omega^f_\psi$ satisfies a third order Picard-Fuchs equation that can be explicitly constructed by the Griffiths-Dwork method. In fact it is more convenient to write down this Picard-Fuchs equation with respect to the parameter $t$ \ref{eq:tpsirelation} instead of $\psi$. From the paper \cite{Nagura}, the Picard-Fuchs equation of $\Omega^f_\psi$ is given by
\begin{equation} \label{eq:k3pfoperator}
\mathcal{D}_3 \Omega^f_\psi=0;~\mathcal{D}_3= \vartheta^3-t(\vartheta+\frac{1}{4})(\vartheta+\frac{1}{2})(\vartheta+\frac{3}{4}),~\vartheta=t\,\frac{d}{dt}.
\end{equation}
Furthermore, the Picard-Fuchs operator $\mathcal{D}_3$ is the symmetric square of a second order linear differential operator $\mathcal{D}_2$ \cite{Hartmann,Nagura}
\begin{equation} \label{eq:k3sqpf}
\mathcal{D}_2=\vartheta^2-t(\vartheta+\frac{1}{8})(\vartheta+\frac{3}{8}).
\end{equation}
Here symmetric square means that if $\pi_0(t)$ and $\pi_1(t)$ are two linearly independent solutions of the operator $\mathcal{D}_2$, then $\pi^2_0(t)$, $\pi_0(t)\pi_1(t)$ and $\pi^2_1(t)$ are three linearly independent solutions of the operator $\mathcal{D}_3$ \cite{Hartmann, Nagura}. 

The operator $\mathcal{D}_2$ \ref{eq:k3sqpf} has three regular singularities at the points
\begin{equation} \label{eq:singularityt}
t=0,1,\infty;
\end{equation}
which are also all the singularities of its symmetric square $\mathcal{D}_3$ \ref{eq:k3pfoperator}. It should be noticed that the fiber of the Fermat pencil \ref{eq:fermatpencilK3} over $t=\infty$ ($\psi=0$) is smooth. In fact, the interesting behavior of this `fake' singularity will be revealed when we study the connections between the Fermat pencil of K3 surfaces and the Legendre family of elliptic curves in Section \ref{sec:mirrormapLegendre}.

The independent solutions of the Picard-Fuchs operator $\mathcal{D}_3$ \ref{eq:k3pfoperator} have been explicitly found in the paper \cite{Nagura}
\begin{equation} \label{eq:d3independentperiods}
\begin{aligned}
W_0(\psi)&=\sum_{n=0}^\infty \frac{(4n)!}{(n!)^4(4 \psi)^{4n}},\\
W_1(\psi)&=\frac{1}{2 \pi i}\left(-4W_0 \cdot \log(4 \psi)+4 \sum_{n=0}^\infty \frac{(4n)!}{(n!)^4(4\psi)^{4n}}[\Psi(4n+1)-\Psi(n+1)] \right),\\
W_2(\psi)&=\frac{1}{(2 \pi i)^2}\Big[ 4^2W_0[\log(4\psi)]^2-2\cdot 4^2 \sum_{n=0}^\infty \frac{(4n)!}{(n!)^4(4\psi)^{4n}}[\Psi(4n+1)-\Psi(n+1)] \cdot \log(4\psi)\\
&\,\,\,\,\,+4^2 \sum_{n=0}^\infty \frac{(4n)!}{(n!)^4(4\psi)^{4n}} \left\lbrace [\Psi(4n+1)-\Psi(n+1)]^2+\Psi'(4n+1)-\frac{1}{4} \Psi'(n+1) \right\rbrace \Big],
\end{aligned}
\end{equation}
which converge in a neighborhood of $\psi=\infty$. Here $\Psi(z)$ is the polygamma function
\begin{equation}
\Psi(z)=\frac{d}{dz} \log \Gamma(z).
\end{equation}
On the other hand, one solution of the operator $\mathcal{D}_2$ \ref{eq:k3sqpf} is given by the hypergeometric function
\begin{equation} \label{eq:fdperiodst}
\pi_0(t)= ~_2F_1(\frac{1}{8},\frac{3}{8};1;t)= 1 +\frac{3 }{64}\,t+\frac{297 }{16384}\,t^2 +\frac{10659 }{1048576}\,t^3+ \cdots.\\
\end{equation}
Moreover, its square is in fact the solution $W_1(\psi)$ in formula \ref{eq:d3independentperiods}, i.e.
\begin{equation}
W_0(\psi)=\pi^2_0(t),~t=\psi^{-4},
\end{equation}
which follows directly from the property that $\mathcal{D}_3$ is the symmetric square of $\mathcal{D}_2$.

\subsection{The construction of the Dwork family}

The mirror family of quartic K3 surfaces is a pencil of K3 surfaces called the Dwork family, which can be constructed from the Fermat pencil \ref{eq:fermatpencilK3} by the Greene-Plesser construction. More explicitly, the abelian group $G$
\begin{equation}
G=\{(a_0,a_1,a_2,a_3)|a_i^4=1,~a_0a_1a_2a_3=1 \}/ \{(a,a,a,a)|a^4=1 \}
\end{equation}
acts freely on the fiber $\mathscr{F}_\psi$ \ref{eq:fermatpencilK3} through
\begin{equation}
(a_0,a_1,a_2,a_3).(X_0,X_1,X_2,X_3)=(a_0X_0,a_1X_1,a_2X_2,a_3X_3).
\end{equation}
Moreover, $G$ is isomorphic to $(\mathbb{Z}/4\mathbb{Z})^2$, and it permutes the 16 singular points of $\mathscr{F}_{\zeta_4^n}$. For $\psi \notin \Sigma$, the quotient variety $\mathscr{F}_\psi/G$ has 6 singularities of type $A_3$. While if $\psi^4=1$, there is an additional singular point of type $A_1$, which is just the quotient of the 16 singular points of $\mathscr{F}_{\zeta_4^n}$ by $G$. If $\psi=\infty$, the quotient $\mathscr{F}_\infty/G$ is a union of hyperplanes, which is isomorphic to $\mathscr{F}_\infty$ \cite{Hartmann, Nagura}.

There exists a minimal simultaneous resolution of the $A_3$ singularities of $\mathscr{F}_\psi, \psi \neq \infty$, after which we obtain a pencil of K3 surfaces called the Dwork family
\begin{equation} \label{eq:mirrorfamilyofK3}
\pi:\mathscr{X} \rightarrow \mathbb{P}^1,
\end{equation}
which is also defined over $\mathbb{Q}$. The details of this mirror construction are left to the papers \cite{Hartmann, Nagura}. The singular fibers of the Dwork family \ref{eq:mirrorfamilyofK3} are also over the points in $\Sigma$, and the singularity of $\mathscr{X}_{\zeta_4^n}$ is a single point of type $A_1$. The Picard number $\rho(\mathscr{X}_\psi) $ of a smooth fiber $\mathscr{X}_\psi$ is $\geq 19$, and a general smooth fiber has Picard number 19 \cite{Hartmann,Huybrechts, Yui}.

The holomorphic twoform $\Omega^f_\psi$ \ref{eq:quartictwoform} is invariant under the action of $G$, hence it defines a nowhere vanishing twoform on the smooth locus of the quotient $\mathscr{F}_\psi/G$. After resolution of singularities, this twoform extends to a nowhere vanishing holomorphic twoform $\Omega_\psi$ on $\mathscr{X}_\psi$, which satisfies the same Picard-Fuchs equation as $\Omega^f_\psi$, i.e. \cite{Nagura}
\begin{equation}
\mathcal{D}_3(\Omega_\psi)=0.
\end{equation}

\subsection{The mirror map}

Recall from \textbf{Remark} \ref{remarkDFK}, $X$ also means the underlying differential manifold structure of a smooth fiber of the Dwork family \ref{eq:mirrorfamilyofK3t}. Since $\Omega_\psi$ satisfies the same Picard-Fuchs equation as $\Omega^f_\psi$, the three independent solutions in formula \ref{eq:d3independentperiods} are the three independent periods of $\Omega_\psi$. From the papers \cite{Hartmann, Nagura}, there exist two integral homology cycles $\Gamma_0,\Gamma_1 \in H_2(X,\mathbb{Z})$ such that
\begin{equation}
W_0(\psi)=\frac{l}{l(2 \pi i)^2}\int_{\Gamma_0} \Omega_\psi,~W_1(\psi)=\frac{1}{l(2 \pi i)^2}\int_{\Gamma_1} \Omega_\psi;~l\in \mathbb{Q}^\times.
\end{equation}
The mirror map $\tau$ for the Dwork family \ref{eq:mirrorfamilyofK3} is given by \cite{Hartmann,Nagura}
\begin{equation} \label{eq:mirrormapDwork}
\tau= \frac{W_1(\psi)}{W_0(\psi)}.
\end{equation}

The values of $\tau$ lie in the upper half plane $\mathbb{H}$. Moreover, with a suitable choice of branch cuts, the special values of $\tau$ at $\psi=0,1,\infty$ are given by \cite{Nagura}
\begin{equation} \label{eq:singularpointscoords}
\tau: 0 \mapsto \frac{-1+i}{2}; ~1 \mapsto \frac{i}{\sqrt{2}}; ~\infty \mapsto \infty.
\end{equation}
So we can say a fundamental domain for $\tau$ is the hyperbolic triangle with vertices $\frac{-1+i}{2}$, $\frac{i}{\sqrt{2}}$ and $\infty$ \cite{Hartmann,Nagura}. The properties of this mirror map and its connections to the $j$-function have also been studied in the paper \cite{Yau}. Later in this paper, we will show that this mirror map is the same as the period map of the Legendre family of elliptic curves. But first let us review the theories of elliptic curves that will be needed in later sections.

\section{An overview of elliptic curves and the Legendre family} \label{sec:legendrefamilyofellipticcurves}

In this section, we will review some elementary properties about elliptic curves. We will also discuss the Legendre family of elliptic curves and the modular lambda function \cite{CarlsonMullerPeters, Chandra}. This section is included here purely to familiar the readers with the notations in later sections.

\subsection{An overview of elliptic curves defined over $\mathbb{Q}$}

First, let us look at the elliptic curves defined over $\mathbb{Q}$. Given an elliptic curve $\mathcal{E}$ defined over $\mathbb{Q}$, it always has an integral model of the form \cite{Silverman,Silverman1}
\begin{equation} \label{eq:weierstrassequation}
\mathcal{E}: y^2+a_1 x y+a_3 y=x^3+a_2 x^2+a_4 x+a_6,~\text{with}~a_1,\cdots, a_6 \in \mathbb{Z},
\end{equation}
which is called integral Weierstrass equation. The discriminant $\Delta$ of this Weierstrass equation \ref{eq:weierstrassequation} is by definition 
\begin{equation}
\Delta=-b_2^2b_8-8b_4^3-27b_6^2+9b_2b_4b_6,
\end{equation}
where $b_i$ is given in terms of $a_i$
\begin{equation}
\begin{aligned}
b_2&=a_1^2+4a_2,\\
b_4&=2a_4+a_1a_3,\\
b_6&=a_3^2+4a_6,\\
b_8&=a_1^2a_6+4a_2a_6-a_1a_3a_4+a_2a_3^2-a_4^2.
\end{aligned}
\end{equation}
The elliptic curve $\mathcal{E}$ \ref{eq:weierstrassequation} is smooth if and only if $\Delta \neq 0$. An elliptic curve $\mathcal{E}$ can have many different integral Weierstrass equations, and a minimal Weierstrass equation is one for which the absolute value $|\Delta|$ is minimal among all Weierstrass models for $\mathcal{E}$. In fact, given an elliptic curve defined over $\mathbb{Q}$, it always has a minimal Weierstrass equation. The $j$-invariant of $\mathcal{E}$ \ref{eq:weierstrassequation} is by definition
\begin{equation}
j(\mathcal{E})=\frac{c_4^3}{\Delta},
\end{equation}
where $c_4$ is defined by
\begin{equation}
c_4=b_2^2-24b_4.
\end{equation}

The endomorphism ring of $\mathcal{E}$, denoted by $\text{End}(\mathcal{E})$, is the ring that consists of all the endomorphisms of $\mathcal{E}$, including those defined over extensions of the base field $\mathbb{Q}$. An elliptic curve $\mathcal{E}$ does not admit complex multiplication (CM) if $\text{End}(\mathcal{E})$ is isomorphic to $\mathbb{Z}$. While $\mathcal{E}$ is said to admit CM if $\text{End}(\mathcal{E})$ is bigger than $\mathbb{Z}$, in which case $\text{End}(\mathcal{E})$ is an order in an imaginary quadratic field. Recall that an \textbf{order} of an algebraic number field $K$ is a sub-ring $\mathcal{O}$ of $\mathcal{O}_K$, the ring of integers of $K$, that is also a $\mathbb{Z}$-module of rank $[K:\mathbb{Q}]$ \cite{LMFDB, Neukirch}. For an elliptic curve defined over $\mathbb{Q}$ that admits CM, its order is one of the 13 orders of class number one \cite{LMFDB,Silverman}. In fact, the property of admitting CM only depends on the $j$-invariants of elliptic curves. For a CM elliptic curve defined over $\mathbb{Q}$, its $j$-invariant is one of the following 13 CM $j$-invariants \cite{LMFDB,Silverman}
\begin{equation}
\begin{aligned}
j=&-262537412640768000, -147197952000, -884736000, -12288000, -884736,\\
& -32768, -3375, 0, 1728, 8000, 54000, 287496, 16581375.
\end{aligned}
\end{equation}

Modulo a prime number $p$, the integral Weierstrass equation \ref{eq:weierstrassequation} defines a curve $\mathcal{E}/\mathbb{F}_p$ over the finite field $\mathbb{F}_p=\mathbb{Z}/p\mathbb{Z}$. $\mathcal{E}$ is said to have good (resp. bad) reduction at $p$ if $\mathcal{E}/\mathbb{F}_p$ is smooth (resp. singular). We will call $p$ a good (resp. bad) prime of $\mathcal{E}$ if $\mathcal{E}$ has good (resp. bad) reduction at $p$. The conductor of $\mathcal{E}$, denoted by $N(\mathcal{E})$, is determined solely by its bad primes, whose precise definition is left to \cite{LMFDB,Silverman,Silverman1}. Let us denote the number of points of $\mathcal{E}/\mathbb{F}_p$ for a good prime $p$ by $\#(\mathcal{E}/\mathbb{F}_p)$, and let $a_p(\mathcal{E})$ be
\begin{equation}
a_p(\mathcal{E})=1+p-\#(\mathcal{E}/\mathbb{F}_p).
\end{equation}
Then the zeta function of $\mathcal{E}$ for the good prime $p$ is of the form \cite{Diamond}
\begin{equation}
\zeta(\mathcal{E},p,T)=\frac{1-a_p(\mathcal{E})T+pT^2}{(1-T)(1-pT)}.
\end{equation}
Appendix \ref{sec:Weilconjectures} contains a short review about zeta functions and Weil conjectures. On the other hand, the \'etale cohomology group $H^1_{\text{\'et}}(\mathcal{E},\mathbb{Q}_\ell)$ is a two dimensional continuous representation of the absolute Galois group $\text{Gal}(\overline{\mathbb{Q}}/\mathbb{Q})$. At a good prime $p$, $H^1_{\text{\'et}}(\mathcal{E},\mathbb{Q}_\ell)$ is unramified, and the characteristic polynomial of the geometric Frobenius is \cite{Diamond,Neukirch}
\begin{equation}
1-a_p(\mathcal{E})T+pT^2=(1-\pi^1_p(\mathcal{E})T)(1-\pi^2_p(\mathcal{E})T),
\end{equation}
where the absolute value of $\pi^i_p(\mathcal{E})$ is $p^{1/2}$. See Appendix \ref{sec:Weilconjectures} for more details. The modularity theorem of elliptic curves tells us that $a_p(\mathcal{E})$ is the $p$-th coefficient of the $\textbf{q}$-expansion of a weight-2 newform with level $N(\mathcal{E})$ \cite{Diamond}. We now give an example that will be important in this paper.
\begin{example} \label{ellipticcurveexample}
The elliptic curve 
\begin{equation}
\mathcal{E}_1:y^2=x^3-x
\end{equation}
is labeled as \textbf{32.a3} in LMFDB. Its $j$-invariant is 1728, and its endomorphism ring is $\mathbb{Z}[-1]$, so it admits CM. The weight-2 newform associated to $\mathcal{E}_1$ is labeled as \textbf{32.2.a.a} in LMFDB.
\end{example}
Given two elliptic curves defined over $\mathbb{Q}$, if they have the same $j$-invariant, then they are isomorphic over a number field \cite{Silverman1}. In fact, their difference is a twist, and an introduction to the theory of twisting can be found in the book \cite{Silverman1}. For example, the $j$-invariant of the following elliptic curve
\begin{equation}
\mathcal{E}_2:y^2=x^3-4x,
\end{equation}
is also 1728. But $\mathcal{E}_2$ is not isomorphic to $\mathcal{E}_1$ over $\mathbb{Q}$, instead they are isomorphic over the quadratic field $\mathbb{Q}(\sqrt{2})$. The weight-2 newform associated to $\mathcal{E}_2$ is labeled as \textbf{64.2.a.a} in LMFDB. The difference between the two modular forms is a twist by the Dirichlet character $(2/\cdot)$, i.e. at a good prime $p$ we have \cite{LMFDB}
\begin{equation}
a_p(\mathcal{E}_1)=a_p(\mathcal{E}_2) \left( \frac{2}{p}\right).
\end{equation}
Here $(2/p)$ is the Legendre symbol \cite{Diamond}.

\subsection{The Legendre family of elliptic curves} \label{sec:legendrefamily}

The Legendre family of elliptic curves has played a very important role in the development of modern mathematics, which is defined by the equation
\begin{equation} \label{eq:legendrefamilyequation}
\mathscr{E}_{\lambda}: y^2=x(x-1)(x-\lambda).
\end{equation}
In a more formal language, the formula \ref{eq:legendrefamilyequation} defines a family of elliptic curves over $\mathbb{P}^1$
\begin{equation} \label{eq:fibrationEP1c}
\pi^e:\mathscr{E} \rightarrow \mathbb{P}^1,
\end{equation}
whose fiber over $\lambda \in \mathbb{P}^1$ is $\mathscr{E}_\lambda$ \ref{eq:legendrefamilyequation}. The singular fibers of the Legendre family are over the points 0, 1 and $\infty$. The $j$-invariant of a smooth fiber $\mathscr{E}_\lambda$ is \cite{Chandra}
\begin{equation}
j(\mathscr{E}_\lambda)=256 \frac{(1-\lambda+\lambda^2)^3}{\lambda^2(1-\lambda)^2}.
\end{equation}

Let us now recall the geometric construction of $\mathscr{E}_\lambda$ from cutting and gluing, and we will closely follow the book \cite{CarlsonMullerPeters}. First, cut the complex plane along the line from 0 to $\lambda$ and the line from 1 to $\infty$. Next take a second copy of the complex plane and cut it along the same lines. Then glue the two copies of complex plane together along the branch cuts. What we have obtained is a torus with a complex structure parameterized by $\lambda$. The readers can consult the book \cite{CarlsonMullerPeters} for more details and pictures. 

The canonical bundle of $\mathscr{E}_{\lambda}$ is trivial, and there exists a nowhere vanishing oneform
\begin{equation}
\omega_\lambda=dx/(2 y). 
\end{equation}
The periods of $\omega_\lambda$ are well-known since the nineteenth century, but their computations are still included here. The underlying differential manifold of a smooth fiber $\mathscr{E}_\lambda$ \ref{eq:legendrefamilyequation} is the torus $T=S^1 \times S^1$. Let us now construct a basis $\{\gamma_0, \gamma_1 \}$ for the homology group $H_1(T,\mathbb{Z}) \simeq \mathbb{Z}^2$ from the branch-cut construction of $\mathscr{E}_\lambda$ in the previous paragraph. Let $\gamma_0$ be the cycle that encircles the line $(1,\infty)$ in one copy of the complex plane $\mathbb{C}$, while let $\gamma_1$ be the circle that is the composite of the line from $1$ to $\lambda$ in the first copy and the line from $\lambda$ to $1$ in the second copy. The dual of $\{\gamma_0, \gamma_1 \}$, denoted by $\{\gamma^0,\gamma^1 \}$, forms a basis of the cohomology group $H^1(T, \mathbb{Z}) \simeq \mathbb{Z}^2$. The integration of the oneform $\omega_\lambda$ over the cycles $\{\gamma_0, \gamma_1 \}$ defines two periods of $\mathscr{E}_\lambda$
\begin{equation}
\int_{\gamma_0} \omega_\lambda=2 \pi \varpi_0(\lambda),~\int_{\gamma_1} \omega_\lambda=2 \pi \varpi_1(\lambda).
\end{equation} 
More explicitly, the two periods $\{\varpi_0(\lambda), \varpi_1(\lambda) \}$ are given by the integrals
\begin{equation} \label{eq:periodsintegration}
\begin{aligned}
\varpi_0(\lambda)&=\frac{1}{\pi}\int_{1}^{\infty} \frac{dx}{\sqrt{x(x-1)(x-\lambda)}},\\
\varpi_1(\lambda)&=\frac{1}{\pi}\int_{1}^{\lambda}\, \frac{dx}{\sqrt{x(x-1)(x-\lambda)}}.
\end{aligned}
\end{equation}
After a change of variable by $x=1/z$, the first integral in the formula \ref{eq:periodsintegration} becomes
\begin{equation}
\varpi_0(\lambda)=\frac{1}{\pi}\int_0^1 \frac{dz}{\sqrt{z(1-z)(1-\lambda z)}}.
\end{equation}
If $\lambda$ lies in a small neighborhood of $0$, we can take a series expansion of the factor $(1-\lambda z)^{-1/2}$, then this integral can be computed order by order. The result is a series expansion of $\varpi_0(\lambda)$
\begin{equation}
\varpi_0(\lambda)=1+\frac{1}{4}\,\lambda + \frac{9}{64}\, \lambda^2+ \cdots.
\end{equation}
The second integral in the formula \ref{eq:periodsintegration} can be evaluated similarly, and it admits an expansion with leading terms 
\begin{equation}
\varpi_1(\lambda)=-\frac{1}{\pi i} \int_{\lambda}^1 \frac{dx}{\sqrt{x(1-x)(x-\lambda)}}=-\frac{1}{\pi i}(4 \log 2- \log \lambda)+ \cdots,
\end{equation}
where the limit of the terms in $\cdots$ is zero when $\lambda \rightarrow 0$. By monodromy consideration, we deduce that $\varpi_1(\lambda)$ must be of the form
\begin{equation}
\varpi_1(\lambda)=\frac{1}{\pi i} (\varpi_0(\lambda) \log \lambda+h(\lambda))-\frac{\log 16}{\pi i} \varpi_0(\lambda),
\end{equation}
where $h(\lambda)$ admits a series expansion in a small neighborhood of $\lambda=0$
\begin{equation}
h(\lambda)=\frac{1 }{2}\,\lambda+\frac{21 }{64}\,\lambda^2+\frac{185 }{768}\,\lambda^3+ \cdots.
\end{equation}

\subsection{The modular lambda function} \label{sec:lambdafnidentities}

The nowhere vanishing holomorphic oneform $\omega_\lambda$ satisfies a well-known second order Picard-Fuchs equation   \cite{CarlsonMullerPeters}
\begin{equation} \label{eq:periodsPFequation}
\lambda (1-\lambda)\frac{d^2\omega_\lambda}{d\lambda^2}+(1-2\lambda) \frac{d\omega_\lambda}{d\lambda}-\frac{\lambda}{4} \,\omega_\lambda =0.
\end{equation}
Hence $\varpi_0(\lambda)$ is given by the hypergeometric function \cite{Chandra,YangSW}
\begin{equation} 
\varpi_0(\lambda)=~ _2F_1(\frac{1}{2},\frac{1}{2};1;\lambda).
\end{equation}
The period $\tau$ of the elliptic curve $\mathscr{E}_\lambda$ is by definition given by the quotient
\begin{equation} \label{eq:periodLegendrefamily}
\tau = \frac{\varpi_1(\lambda) }{\varpi_0(\lambda) }.
\end{equation}
In a small neighborhood of $\lambda=0$, $\tau$ is of the form
\begin{equation} \label{eq:defnoftau}
\tau=\frac{1}{\pi i}\, \left(\log \lambda+\frac{h(\lambda) }{\varpi_0(\lambda) } \right)-\frac{\log 16}{\pi i}.
\end{equation}
The underlying complex torus of the elliptic curve $\mathscr{E}_{\lambda}$ is isomorphic to the quotient of $\mathbb{C}$ by the rank-2 lattice generated by $1$ and $\tau$ \cite{CarlsonMullerPeters}. Let $\lambda$ be the coordinate of $\mathbb{P}^1$, then $\tau$ defines a map \cite{Chandra}
\begin{equation}
\tau: \mathbb{P}^1 \rightarrow \mathbb{H} \cup \{\infty\},
\end{equation}
which is called the period map of the Legendre family. The inverse of $\tau$ is the famous modular lambda function, which generates the function field of the modular curve $X(2)$, i.e. it is a Hauptmodul for $X(2)$ \cite{Chandra,Diamond}. In this paper, we will let $q$ be
\begin{equation}
q:=\exp \pi i \tau.
\end{equation}
Formula \ref{eq:defnoftau} implies
\begin{equation} \label{eq:qintermsoflambda}
q=\frac{1}{16} \,\lambda\, \exp \left(h(\lambda)/\varpi_0(\lambda)\right),
\end{equation}
and it can be inverted order by order which gives us the series expansion of $\lambda$ with respective to $q$ \cite{Chandra, YangSW}
\begin{equation} \label{eq:modularLambdafunction}
\lambda(\tau)=16q-128 q^2+704q^3-3072 q^4+11488q^5-38400q^6+ \cdots.
\end{equation}

Furthermore, it is well-known that the period $\varpi_0(\lambda)$ can also be expressed in terms of the theta function $\theta_3$ \cite{Chandra, YangSW}
\begin{equation} \label{eq:varpi0theta3}
\varpi_0(\lambda)=\theta_3^2(0,q).
\end{equation}
From this identity, we also have \cite{Chandra,YangSW}
\begin{equation} \label{eq:theta2and4varpi_0}
\theta_2^4(0,q)=\lambda \varpi^2_0(\lambda),~\theta_4^4(0,q)=(1-\lambda) \varpi^2_0(\lambda).
\end{equation}
From \cite{YangSW}, there is another interesting identity involving $\lambda(\tau)$ and $\varpi_0(\lambda)$ of the form
\begin{equation} \label{eq:differentialsvarpi0}
\frac{1}{\pi i} \frac{d \lambda}{d \tau}=\lambda(1-\lambda)\varpi_0^2(\lambda).
\end{equation}
With these identities at hands, we are ready to study the connections between the mirror map \ref{eq:mirrormapDwork} of the Dwork family and the period map \ref{eq:periodLegendrefamily} of the Legendre family.

\section{The mirror map of the Dwork family and periods of the Legendre family} \label{sec:mirrormapLegendre}

In this section, we will study the connections between the mirror map of the Dwork family and the periods of the Legendre family. More concretely, we will explicitly express the solutions $W_i(\psi)$ \ref{eq:d3independentperiods} in terms of the periods $\varpi_i(\lambda)$ of the Legendre family. Then we will show the mirror map \ref{eq:mirrormapDwork} of the Dwork family is the same as the period map \ref{eq:periodLegendrefamily} of the Legendre family. Based on this result, we will discuss the modularities of the counting functions for K3 surfaces from the mirror symmetry point of view, which shed further lights on this subject. The crucial tools in this section are the quadratic transformations of hypergeometric functions \cite{HG1,HG2,HG3}.

\subsection{The quadratic transformations of periods} \label{sec:quadraperiods}

First, we will need the following two quadratic transformations of hypergeometric functions \cite{HG1,HG2,HG3}
\begin{equation}
\begin{aligned}
_2F_1(\frac{1}{2},\frac{1}{2};1;z)&=(1-z)^{-1/4} \,_2F_1(\frac{1}{4},\frac{1}{4};1;\frac{z^2}{4z-4}),\\
_2F_1(\frac{1}{4},\frac{1}{4};1;z)&=(1-z)^{-1/4}\,_2F_1(\frac{1}{8},\frac{3}{8};1;-\frac{4z}{(1-z)^2}),
\end{aligned}
\end{equation}
which are over the regions when both sides are well-defined. We will however mainly focus on a smooth neighborhood of $z=0$. The composition of these two quadratic transformations gives us 
\begin{equation} \label{eq:quadrticTRHyper}
_2F_1(\frac{1}{2},\frac{1}{2};1;z)=(1-\frac{z}{2})^{-1/2}\,_2F_1(\frac{1}{8},\frac{3}{8};1;-\frac{16(z-1)z^2}{(z-2)^4}).
\end{equation}
More concretely, the power series expansions of the two sides of \ref{eq:quadrticTRHyper} in a small neighborhood of $z=0$ are the same. Now let us define a transformation between the variables $t$ ($=\psi^{-4}$) and $\lambda$ by the following algebraic equation
\begin{equation} \label{eq:quadratictransformation}
t=\lambda^2(1-\lambda)\left(1-\frac{\lambda}{2} \right)^{-4}.
\end{equation}
This equation \ref{eq:quadratictransformation} defines a map from $\mathbb{P}^1$ (with coordinate $\lambda$) to $\mathbb{P}^1$ (with coordinate $t$)
\begin{equation} \label{eq:qtmap}
t:\mathbb{P}^1 \rightarrow \mathbb{P}^1,
\end{equation}
which is a ramified covering map with degree 4. The three singular points $t=0,1,\infty$ of $\mathcal{D}_3$ correspond to
\begin{equation} \label{eq:singularitytransformations}
\begin{aligned}
t=0 & \iff \lambda=0,1, \infty;\\
t=1 & \iff \lambda=\pm 2 \sqrt{2}-2;\\
t= \infty & \iff \lambda=2.
\end{aligned}
\end{equation}
The map $t$ \ref{eq:qtmap} has four ramification points: $\lambda=0, \pm 2 \sqrt{2}-2$ and 2, where the ramification index of $0, \pm 2 \sqrt{2}-2$ is 2 and that of $2$ is $4$.

The fiber $\mathscr{X}_\psi$ of the Dwork family \ref{eq:mirrorfamilyofK3} is isomorphic to the fiber $\mathscr{X}_{\zeta_4\psi}$, hence the Dwork family \ref{eq:mirrorfamilyofK3} descends to a family over $\mathbb{P}^1$ with parameter $t$
\begin{equation} \label{eq:mirrorfamilyofK3t}
\pi^t:\mathscr{X} \rightarrow \mathbb{P}^1.
\end{equation}
The pull-back of this family \ref{eq:mirrorfamilyofK3t} along the map $t$ \ref{eq:qtmap} gives us a commutative diagram
\begin{equation} \label{eq:pullbackCD}
\begin{tikzcd}
\widetilde{\mathscr{X}} \arrow[r] \arrow[d,"\tilde{\pi}"] & \mathscr{X} \arrow[d,"\pi^t"] \\
\mathbb{P}^1 \arrow[r,"t"] & \mathbb{P}^1
\end{tikzcd}.
\end{equation}
The new family $\tilde{\pi}$ in this commutative diagram is also a pencil of K3 surfaces over $\mathbb{P}^1$
\begin{equation} \label{eq:pullbackk3family}
\tilde{\pi}:\widetilde{\mathscr{X}} \rightarrow \mathbb{P}^1,
\end{equation}
which will be crucial in this paper. Later we will show that in a sense this family is a `more suitable' mirror family for quartic K3 surfaces \ref{eq:kahlerside}.

\subsection{The mirror map is the period map} \label{sec:mirrormapisperiodmap}

Intuitively, we can pull everything on the family \ref{eq:mirrorfamilyofK3t} back to the family \ref{eq:pullbackk3family}. For example, up to an overall factor $16 \lambda$, the operator $\mathcal{D}_2$ \ref{eq:k3sqpf} pulls back to a second order differential operator $\widetilde{\mathcal{D}}_2$
\begin{equation} \label{eq:d2operatorlambda}
\widetilde{\mathcal{D}}_2=\lambda (1-\lambda)(2-\lambda)^2 \frac{d^2}{d\lambda^2}+(2-\lambda)(2-4\lambda+\lambda^2) \frac{d}{d\lambda}-\frac{3}{4} \lambda,
\end{equation}
which has regular singularities at the points
\begin{equation} \label{eq:singularitylambda}
\lambda=0,1,2,\infty.
\end{equation}
Under the map $t$ \ref{eq:qtmap}, the solution $\pi_0(t)$ \ref{eq:fdperiodst} pulls back to
\begin{equation}
\pi_0^K(\lambda)=(1-\frac{\lambda}{2})^{1/2} \,_2F_1(\frac{1}{2},\frac{1}{2};1;\lambda)=(1-\frac{\lambda}{2})^{1/2} \, \varpi_0(\lambda),
\end{equation}
where we have used the identity \ref{eq:quadrticTRHyper}. Then from Section \ref{sec:legendrefamily}, we learn that a second independent solution of $\widetilde{\mathcal{D}}_2$ \ref{eq:d2operatorlambda} is given by
\begin{equation}
\pi_1^K(\lambda)=(1-\frac{\lambda}{2})^{1/2} \varpi_1(\lambda).
\end{equation}
The pull-back of the operator $\mathcal{D}_3$ \ref{eq:k3pfoperator}, denoted by $\widetilde{\mathcal{D}}_3$, is the symmetric square of $\widetilde{\mathcal{D}}_2$ \ref{eq:d2operatorlambda} (up to an overall factor). The holomorphic twoform $\Omega_\psi$ on $\mathscr{X}_\psi$ induces a holomorphic twoform $\widetilde{\Omega}_\lambda$ on $\widetilde{\mathscr{X}}_\lambda$. Independent solutions $\{\Pi_0(\lambda), \Pi_1(\lambda), \Pi_2(\lambda) \}$ of $\widetilde{\mathcal{D}}_3$ are given by
\begin{equation} \label{eq:canonicalperiodsK3}
\begin{aligned}
\Pi_0(\lambda)&=(\pi^K_0(\lambda))^2=(1-\frac{\lambda}{2})\, \varpi^2_0(\lambda), \\
\Pi_1(\lambda) &=\pi^K_0(\lambda) \pi^K_1(\lambda)=(1-\frac{\lambda}{2})\, \varpi_0(\lambda)\varpi_1(\lambda) , \\
\Pi_2(\lambda)&=(\pi^K_1(\lambda))^2=(1-\frac{\lambda}{2})\, \varpi^2_1(\lambda).
\end{aligned}
\end{equation}
\begin{remark}
In this paper, we have assumed that a suitable branch cut has been chosen for the transformation defined by the equation \ref{eq:quadratictransformation}. In this section, we have focused on a small neighborhood of $t=0$ and $\lambda=0$, where we have used the following expansion of equation \ref{eq:quadratictransformation}
\begin{equation}
t=\lambda^2+O(\lambda^3).
\end{equation}
\end{remark}
The crucial observation is that under the transformation \ref{eq:quadratictransformation} we have
\begin{equation}
W_0(\psi)=\Pi_0(\lambda),~W_1(\psi)=\Pi_1(\lambda),
\end{equation}
which can be obtained from the limit behaviors of $W_i(\psi)$ and $\varpi_i(\lambda)$ for
\begin{equation}
\lambda \rightarrow 0 ~\text{and}~\psi \rightarrow \infty.
\end{equation}
Therefore we immediately obtain a crucial property about the mirror map \ref{eq:mirrormapDwork} of the Dwork family
\begin{equation}\label{eq:newmirrormap}
\tau=\frac{W_1(\psi)}{W_0(\psi)}=\frac{\Pi_1(\lambda)}{\Pi_0(\lambda)}=\frac{\varpi_1(\lambda)}{\varpi_0(\lambda)}.
\end{equation}
Namely, the mirror map \ref{eq:mirrormapDwork} of the Dwork family is the same as the period map \ref{eq:periodLegendrefamily} of the Legendre family. This will provide a very important link between the mirror symmetry of K3 surfaces and the Legendre family of elliptic curves. The new family of K3 surfaces \ref{eq:pullbackk3family} can also be considered as the mirror family of quartic K3 surfaces \ref{eq:kahlerside}, and it is actually `better' from a number theoretic point of view!

\subsection{The modularities of counting functions for K3 surfaces}

From the mirror symmetry point of view, the results in Section \ref{sec:mirrormapisperiodmap} will provide philosophical interpretations to an interesting phenomenon that counting functions for K3 surfaces are modular. The philosophy of mirror symmetry says that under the mirror map $\tau$ \ref{eq:newmirrormap}, the counting functions for K3 surfaces on the K\"ahler side correspond to rational expressions of
\begin{equation}
\lambda,~\Pi_0(\lambda), ~d \lambda/d \tau,
\end{equation} 
on the complex side. From Section \ref{sec:lambdafnidentities}, the latter is clearly modular!

We now use a famous example to illustrate this point. The counting function of BPS states in IIB string theory for a K3 surface $X$ (times $\mathbb{R} \times S^1$) has been explicitly worked out in the paper \cite{Vafa}, which is given by
\begin{equation} \label{eq:BPSstatescounting}
\textbf{q}^{-1} \sum_n \chi(\text{Hilb}^n(X))\textbf{q}^n=\frac{1}{\eta^{24}(\tau)}=\frac{1}{\Delta};~\textbf{q}:=\exp (2 \pi i \tau).
\end{equation}
Here $\Delta$ is called the Ramanujan tau function, which can also be expressed in terms of theta functions as \cite{Chandra}
\begin{equation} \label{eq:DeltaexpressionTheta}
\Delta=2^{-8} \theta^8_2(0,q)\theta^8_3(0,q)\theta^8_4(0,q).
\end{equation}
\begin{remark}
In this paper, we will use the notation $\textbf{q}$ to mean $\exp (2 \pi i \tau)$, which is differential from $q=\exp(\pi i \tau)$.
\end{remark}
This counting function \ref{eq:BPSstatescounting} has an alternative derivation, which corresponds to the counting of nodal curves in K3 surfaces \cite{Zaslow}
\begin{equation}
\textbf{q}^{-1} \sum_g \chi(\mathcal{M}_g^H)\textbf{q}^n=\frac{1}{\eta^{24}(\tau)}.
\end{equation}
Here $\mathcal{M}_g^H$ is the moduli space that describes a choice of a holomorphic Riemann surface in K3 surface with genus $g$ and a flat $U(1)$ bundle. The interested readers are referred to the paper \cite{Zaslow} for more details.

On the complex side, using the identities in Section \ref{sec:lambdafnidentities} and formula \ref{eq:DeltaexpressionTheta}, $\Delta$ can be expressed as
\begin{equation}
\Delta=\frac{1}{4}\frac{\lambda^2(1-\lambda)^2}{(\lambda-2)^6} \Pi_0^6(\lambda).
\end{equation}
But of course there are other expressions of $\Delta$ in terms of $\lambda$, $\Pi_0(\lambda)$ and $d \lambda /d \tau$. The upshot is that under the mirror map \ref{eq:newmirrormap}, the counting function \ref{eq:BPSstatescounting} corresponds to
\begin{equation}
\textbf{q}^{-1} \sum_n \chi(\text{Hilb}^n(X))\textbf{q}^n=\frac{4 (\lambda-2)^6}{\lambda^2(1-\lambda)^2} \frac{1}{\Pi_0^6(\lambda)},
\end{equation}
the form of which is certainly what the mirror symmetry of K3 surfaces has predicted. Furthermore, it is very interesting to see whether the results in this section can be applied to study the general counting functions for K3 surfaces. 

\section{Connections with Shioda-Inose structures? } \label{sec:ShiodaInose}

In this section, we will explore the potential connections between the results of Section \ref{sec:mirrormapLegendre} and the Shioda-Inose structures of smooth fibers of the Fermat pencil and Dwork family.

Recall from \textbf{Remark} \ref{remarkDFK} that $X$ also means the underlying differential manifold of a K3 surface. From \cite{Hartmann}, there exist integral homology cycles $h,e,f \in H^2(X,\mathbb{Z})$ such that the holomorphic twoform $\widetilde{\Omega}_\lambda$ on the smooth fiber $\widetilde{\mathscr{X}}_\lambda$ of the family \ref{eq:pullbackk3family} admits an expansion
\begin{equation} \label{eq:twoformexpansion}
\widetilde{\Omega}_\lambda=l(2 \pi i)^2 \left( \Pi_1(\lambda)h-\Pi_0(\lambda)e+2\Pi_2(\lambda)f \right),
\end{equation}
where $l$ is a nonzero rational constant. Moreover, the only nontrivial cup-product pairings between $h$, $e$ and $f$ are 
\begin{equation} \label{eq:hefpairing}
\langle h,h \rangle=4,~ \langle e, f\rangle=\langle f,e\rangle=1.
\end{equation}
For simplicity, let the free $\mathbb{Z}$-module generated by $h$, $e$ and $f$ be $L_1$
\begin{equation}
L_1=\mathbb{Z}h \oplus \mathbb{Z}e \oplus \mathbb{Z}f.
\end{equation}
Given a $\lambda$ such that 
\begin{equation}
\lambda \neq 0, 1, \pm 2 \sqrt{2}-2,2, \infty,
\end{equation}
if the Picard number of $\widetilde{\mathscr{X}}_\lambda$ is 19, then the transcendental lattice $T(\widetilde{\mathscr{X}}_\lambda)$ of $\widetilde{\mathscr{X}}_\lambda$ is just $L_1$. While if the Picard number of $\widetilde{\mathscr{X}}_\lambda$ is 20, then the transcendental lattice $T(\widetilde{\mathscr{X}}_\lambda)$ of $\widetilde{\mathscr{X}}_\lambda$ is a rank-2 sub-lattice of $L_1$. 

On the other hand, the direct product of the Legendre family \ref{eq:fibrationEP1c}, i.e. $\mathscr{A}=\mathscr{E} \times \mathscr{E}$, is a family of complex surfaces over $\mathbb{P}^1$, and the fiber $\mathscr{A}_\lambda$ is just the direct product $\mathscr{E}_\lambda \times \mathscr{E}_\lambda$. The underlying differential manifold of a smooth fiber $\mathscr{A}_\lambda$ is the direct product of torus, i.e. $T \times T$. A smooth fiber $\mathscr{A}_\lambda$ has three rationally independent algebraic cycles
\begin{equation} \label{eq:algedivisorsprod}
\mathscr{E}_\lambda \times 0, ~0\times \mathscr{E}_\lambda, ~\Delta_\lambda,
\end{equation}
where $\Delta_\lambda$ is the diagonal of $\mathscr{E}_\lambda \times \mathscr{E}_\lambda$. Therefore the Picard number of $\mathscr{A}_\lambda$ is $\geq 3$. The integral cohomology group $H^2(T \times T,\mathbb{Z})$ is a free $\mathbb{Z}$-module of rank 6 with a unimodular cup-product pairing. Under this pairing, the orthogonal complement of the three algebraic cycles in formula \ref{eq:algedivisorsprod} is the lattice
\begin{equation} \label{eq:abeliansurfacesgenerators}
L_2=\mathbb{Z}(\gamma^0\otimes \gamma^0) \oplus \mathbb{Z}(\gamma^0\otimes \gamma^1+\gamma^1 \otimes \gamma^0) \oplus \mathbb{Z}(\gamma^1 \otimes \gamma^1).
\end{equation}
The only nontrivial pairings between the three generators of $L_2$ are
\begin{equation}
\begin{aligned}
&\langle \gamma^0\otimes \gamma^0,\gamma^1\otimes \gamma^1\rangle=\langle \gamma^1\otimes \gamma^1,\gamma^0\otimes \gamma^0\rangle=1.\\
&\langle \gamma^0\otimes \gamma^1+\gamma^1 \otimes \gamma^0,\gamma^0\otimes \gamma^1+\gamma^1 \otimes \gamma^0 \rangle=-2.
\end{aligned}
\end{equation}
Recall from Section \ref{sec:legendrefamily} that $\{\gamma^0,\gamma^1 \}$ is a basis of $H^1(T,\mathbb{Z})$. The nowhere vanishing holomorphic twoform on $\mathscr{A}_\lambda$ is given by the tensor product $\omega_\lambda \otimes \omega_\lambda$, where $\omega_\lambda$ is the nowhere vanishing holomorphic oneform on the elliptic curve $\mathscr{E}_\lambda$. From Section \ref{sec:legendrefamily}, $\omega_\lambda \otimes \omega_\lambda$ admits an expansion of the form
\begin{equation} \label{eq:abeliansurfaceTwoformexpansion}
\omega_\lambda \otimes \omega_\lambda=(2 \pi)^2\left[ \varpi_0^2(\lambda) \gamma^0\otimes \gamma^0 +\varpi_0(\lambda)\varpi_1(\lambda )(\gamma^0\otimes \gamma^1+\gamma^1 \otimes \gamma^0)+ \varpi_1^2(\lambda) \gamma^1 \otimes \gamma^1 \right].
\end{equation}

The pure Hodge structure on the transcendental lattice $T(\widetilde{\mathscr{X}}_\lambda)$ is determined by the holomorphic twoform $\widetilde{\Omega}_\lambda$, or equivalently its expansion \ref{eq:twoformexpansion}. Similarly, the pure Hodge structure on the transcendental lattice $T(\mathscr{A}_\lambda)$ is determined by the holomorphic twoform $\omega_\lambda \otimes \omega_\lambda$, or equivalently its expansion \ref{eq:abeliansurfaceTwoformexpansion}. From formula \ref{eq:canonicalperiodsK3}, we learn that the pure Hodge structure on $T(\widetilde{\mathscr{X}}_\lambda)$ is isomorphic to that on $T(\mathscr{A}_\lambda)$, but in general it is not a Hodge isometry. Therefore a very interesting question is about the connections between the Shioda-Inose structure of $\widetilde{\mathscr{X}}_\lambda$ and the geometry of the complex surface $\mathscr{A}_\lambda$.

\section{The zeta functions of smooth fibers of the Fermat pencil} \label{sec:zetafunctionsconnections}

In this section, we will look at the potential relations between the zeta functions of $\mathscr{A}_\lambda$, i.e. $\mathscr{E}_\lambda \times \mathscr{E}_\lambda$, and the zeta functions of a rational model of the smooth fiber $\mathscr{F}_{\psi(\lambda)}$ of the Fermat pencil.

\subsection{The pull back of the Fermat pencil}

We can apply the constructions in Section \ref{sec:quadraperiods} to the Fermat pencil and obtain a family of K3 surface with parameter $\lambda$. However, we find it more convenient to write everything down explicitly. More precisely, the equation \ref{eq:quadratictransformation} defines $\psi$ as a multivalued function of $\lambda$
\begin{equation} \label{eq:psilambdaTran}
\psi(\lambda)=\lambda^{-\frac{1}{2}} (1-\lambda)^{-\frac{1}{4}}\left(1-\frac{\lambda}{2} \right),
\end{equation}
where we assume that a suitable branch cut has been chosen. Now the Fermat pencil of K3 surfaces \ref{eq:fermatpencilK3} becomes
\begin{equation} \label{eq:fermatpencilK3Lambda}
\mathscr{F}_{\psi(\lambda)}: \{f_{\psi(\lambda)}=0 \}\subset \mathbb{P}^3, 
\end{equation}
where the quartic polynomial $f_{\psi(\lambda)}$ is
\begin{equation}
f_{\psi(\lambda)}=X_0^4+X_1^4+X_2^4+X_3^4-4\psi(\lambda) \,X_0X_1X_2X_3.
\end{equation}
Similarly, we have a meromorphic threeform $\Theta_{\psi(\lambda)}$ whose residue defines a nowhere vanishing holomorphic twoform $\Omega^f_{\psi(\lambda)}$ on $\mathscr{F}_{\psi(\lambda)}$ that satisfies the same Picard-Fuchs equation as the twoform $\widetilde{\Omega}_\lambda$ on $\widetilde{\mathscr{X}}_\lambda$, i.e.
\begin{equation}
\widetilde{\mathcal{D}}_3 \,\Omega^f_{\psi(\lambda)}=0.
\end{equation}
Moreover, there exist cohomological elements $e_i \in H^2(X,\mathbb{Q})$ such that 
\begin{equation} \label{eq:fermatomegeexpansion}
\Omega^f_{\psi(\lambda)}=l_1(2 \pi i)^2 \left( \Pi_0(\lambda) e_0+\Pi_1(\lambda) e_1+\Pi_2(\lambda) e_2\right),~ l_1 \in \mathbb{Q}^\times.
\end{equation}
\begin{remark} \label{RemarkShiodaInoseFermat}
Similarly from Section \ref{sec:ShiodaInose}, the pure Hodge structure on the transcendental lattice $T(\mathscr{F}_{\psi(\lambda)})$ of a smooth fiber $\mathscr{F}_{\psi(\lambda)}$ of the Fermat pencil is isomorphic to that on $T(\mathscr{A}_\lambda)$.  
\end{remark}

Intuitively, we will say $\mathscr{E}_\lambda$ of the Legendre family is the \textbf{elliptic partner} of the K3 surface $\mathscr{F}_{\psi(\lambda)}$. It is interesting to notice that the special fibers of the Fermat pencil at $\psi =0,1, \infty$ admit very interesting elliptic partners:
\begin{enumerate}
\item When $\psi=0$, we have the famous Fermat quartic
\begin{equation}\label{eq:fermatquartic}
\mathscr{F}_0:\{ X_0^4+X_1^4+X_2^4+X_3^4=0 \} \subset \mathbb{P}^3.
\end{equation}
From formula \ref{eq:psilambdaTran}, $\psi=0$ corresponds to $\lambda=2$, and the smooth fiber of the Legendre family over $\lambda=2$ is 
\begin{equation}
y^2=x(x-1)(x-2),
\end{equation}
whose Weierstrass integral model is just $\mathscr{E}_1$ in \textbf{Example} \ref{ellipticcurveexample}, i.e. \textbf{32.a3} in LMFDB.

\item When $\psi=1$,  we have the singular surface
\begin{equation}
\mathscr{F}_1:\{ X_0^4+X_1^4+X_2^4+X_3^4-4X_0X_1X_2X_3=0 \} \subset \mathbb{P}^3.
\end{equation}
From formula \ref{eq:psilambdaTran}, $\psi=1$ corresponds to $\lambda= 2 \sqrt{2}-2$. The smooth fiber of the Legendre family over $\lambda= 2 \sqrt{2}-2$ is the elliptic curve
\begin{equation}
y^2=x(x-1)(x-( 2 \sqrt{2}-2)),
\end{equation}
both of which are smooth elliptic curves defined over $\mathbb{Q}(\sqrt{2})$ with $j$-invariant 8000.
\item When $\psi=\infty$, we have a union of four complex planes 
\begin{equation}
\mathscr{F}_\infty: \{X_0 X_1X_2X_3=0 \} \subset \mathbb{P}^3.
\end{equation}
From formula \ref{eq:psilambdaTran}, $\psi=\infty$ corresponds to $\lambda=0,1,\infty$, and the fibers of the Legendre family over $\lambda=0,1,\infty$ are just the singular fibers of it.
\end{enumerate}

\subsection{The properties of zeta functions}
For simplicity, let us assume $\lambda \in \mathbb{Q}$ and
\begin{equation}
\lambda \neq 0, 1,2, \infty,
\end{equation}
$\mathscr{A}_\lambda$ is defined over $\mathbb{Q}$. In this section, a rational model for the smooth fiber $\mathscr{F}_{\psi(\lambda)}$ is chosen to be the one given by formula \ref{eq:rationalfpsi}. The transcendental lattice $T(\mathscr{A}_\lambda)$ (resp. $T(\mathscr{F}_{\psi(\lambda)})$) generates a continuous representation of $\text{Gal}(\overline{\mathbb{Q}}/\mathbb{Q})$, which will be denoted by $V^a_\lambda$ (resp. $V^f_{\psi(\lambda)}$). From \textbf{Remark} \ref{RemarkShiodaInoseFermat}, the pure Hodge structure on the rational vector space $T(\mathscr{F}_{\psi(\lambda)})\otimes \mathbb{Q}$ is isomorphic to that on $T(\mathscr{A}_\lambda)\otimes \mathbb{Q}$. Hence from the Hodge conjecture, we learn that there exists a number field $K$ such that $V^a_\lambda$ is isomorphic to $V^f_{\psi(\lambda)}$ as representations of $\text{Gal}(\overline{\mathbb{Q}}/K)$ \cite{Shamit, KimYang}. This property immediately implies that there may exist interesting relations between the zeta functions of $V^a_\lambda$ and that of $V^f_{\psi(\lambda)}$ at good primes. 

Given an elliptic curve $\mathscr{E}_\lambda$ with $\lambda \in \mathbb{Q}$, the zeta function of $H^1_{\text{\'et}}(\mathscr{E}_\lambda,\mathbb{Q}_\ell)$ at a good prime $p$ is a quadratic polynomial
\begin{equation} \label{eq:zetafunctionelambda}
1-a_p(\mathscr{E}_\lambda)T+pT^2=(1-\pi^1_p(\mathscr{E}_\lambda)T)(1-\pi^2_p(\mathscr{E}_\lambda)T).
\end{equation}
See Section \ref{sec:legendrefamilyofellipticcurves} for more details. The symmetric square of formula \ref{eq:zetafunctionelambda} is by definition 
\begin{equation}
(1-\pi^1_p(\mathscr{E}_\lambda)\pi^2_p(\mathscr{E}_\lambda)T)(1-(\pi^1_p(\mathscr{E}_\lambda))^2T)(1-(\pi^2_p(\mathscr{E}_\lambda))^2T),
\end{equation}
which simplifies to
\begin{equation} \label{eq:sqzetafunctionsEC}
(1-pT)(1-(a^2_p(\mathscr{E}_\lambda)-2p)T+p^2T^2).
\end{equation}
This cubic polynomial \ref{eq:sqzetafunctionsEC} is a factor of the zeta function of $H^2_{\text{\'et}}(\mathscr{E}_\lambda \times \mathscr{E}_\lambda,\mathbb{Q}_\ell)$ that corresponds to the lattice $L_2$ \ref{eq:abeliansurfacesgenerators}. So one can ask whether there exists a Dirichlet character $\chi_\lambda$ depending one $\lambda$ such that the twisted cubic polynomial
\begin{equation}
(1-\chi_\lambda(p)pT)(1-\chi_\lambda(p)(a^2_p(\mathscr{E}_\lambda)-2p)T+p^2T^2)
\end{equation}
is a factor of the zeta function of $H^2_{\text{\'et}}(\mathscr{F}_{\psi(\lambda)},\mathbb{Q}_\ell)$  \cite{Ono,Yui2011}? We will not pursue this interesting question further in this paper, while the readers are referred to the paper \cite{Ono} for more details about the computations of zeta functions of a pencil of K3 surfaces using that of elliptic curves. In the rest of this part, we will focus on the case of the Fermat quartic \ref{eq:fermatquartic}.

\begin{remark}
The discussions in this section also apply to the zeta functions of smooth fibers of the Dwork family.
\end{remark}

\section{The Fermat quartic and Deligne's conjecture} \label{sec:DeligneconjectureFermatquartic}

In this section, we will compute the periods of the holomorphic twoform on the Fermat quartic $\mathscr{F}_0$. Then we will discuss the relations between the modularity of the Fermat quartic $\mathscr{F}_0$ and that of the elliptic curve \textbf{32.a3} in LMFDB \cite{LMFDB}
\begin{equation} \label{eq:lmfdb32a3}
y^2=x^3-x.
\end{equation}
We will also apply the method developed in \cite{YangDeligne} to compute Deligne's periods of the Fermat quartic $\mathscr{F}_0$ and (numerically) verify that they satisfy Deligne's conjecture on the special values of $L$-functions at critical integral points \cite{DeligneL}. In this section, we will need the theory of pure motives, which has been briefly reviewed in the papers \cite{Shamit,KimYang}.

\subsection{The periods of the Fermat quartic}

First, let us compute the periods of the Fermat quartic $\mathscr{F}_0$. In the construction of the twoform $\Omega^f_\psi$ on $\mathscr{F}_\psi$ in Section \ref{sec:canonicalperiodsoffermatpencil}, there is an additional factor $\psi$, therefore $\Omega^f_\psi$ becomes 0 on the Fermat quartic $\mathscr{F}_0$. This defect can also be seen from the values of the periods $\Pi_i(\lambda)$ \ref{eq:canonicalperiodsK3} at $\lambda=2$
\begin{equation}
\Pi_0(2)=\Pi_1(2)=\Pi_2(2)=0.
\end{equation}
It can be cured by defining the meromorphic threeform $\Theta_F$ to be
\begin{equation}
\Theta_F=\sum_{i=0}^3 (-1)^i \frac{x_i \,dx_0 \wedge \cdots \wedge \widehat{dx_i} \wedge \cdots \wedge d x_3}{f_0},
\end{equation}
whose residue along $\mathscr{F}_0$ defines a nowhere vanishing holomorphic twoform on the Fermat quartic. More explicitly, take the open subset of $\mathscr{F}_0$ defined by $x_3=1$, then the residue of $\Theta_F$ is equal to \cite{MarkGross}
\begin{equation} \label{eq:meroHoloTwoFermat}
\Omega_F= \frac{dx_0 \wedge dx_1 }{4 x^3_2} \Big|_{\mathscr{F}_0}.
\end{equation}
Similarly, it can be explicitly shown that the meromorphic twoform \ref{eq:meroHoloTwoFermat}, which is a priori only defined on $ x_2 \neq 0$, extends to a global nowhere vanishing twoform on $\mathscr{F}_0$  \cite{MarkGross}. It is very important that $\Omega_F$ is defined over $\mathbb{Q}$, and it spans the algebraic de Rham cohomology group $H^2_{\text{dR}}(\mathscr{F}_0)$ \cite{Shamit,KimYang,YangDeligne}.

The periods of $\Omega_F$ can be found from that of $\Omega^f_\psi$, i.e. $\Pi_i(\lambda)$  \ref{eq:canonicalperiodsK3}. More precisely, the twoform $\Omega_F$ is the limit of $\Omega^f_\psi/\psi$ at $\psi=0$, hence we have
\begin{equation}
\Omega_F=\lim_{\psi \rightarrow 0} \Omega^f_\psi/\psi=\lim_{\lambda \rightarrow 2} \Omega^f_{\psi(\lambda)}/\psi(\lambda).
\end{equation}
Then from formulas \ref{eq:canonicalperiodsK3}, \ref{eq:psilambdaTran} and \ref{eq:fermatomegeexpansion}, we immediately obtain the following crucial expansion of $\Omega_F$
\begin{equation} \label{eq:fermatTwoformExpansion}
\Omega_F=l_1(1+i)(2 \pi i)^2 \left( \varpi^2_0(2) e_0+\varpi_0(2)\varpi_1(2) e_1+\varpi^2_1(2) e_2\right),~ l_1 \in \mathbb{Q}^\times,e_i \in H^2(X,\mathbb{Q}).
\end{equation}
The pure Hodge structure on the transcendental lattice $T(\mathscr{F}_0)$ of the Fermat quartic \ref{eq:fermatquartic} is uniquely determined by the expansion \ref{eq:fermatTwoformExpansion}, therefore it is isomorphic to the pure Hodge structure on the transcendental lattice $T(\mathscr{A}_2)$. Let us now look at the \'etale cohomological counterpart of this property, e.g. zeta functions.

\subsection{The modular form of the Fermat quartic} \label{sec:ModularFermat}

The Fermat quartic $\mathscr{F}_0$ \ref{eq:fermatquartic} is perhaps the earliest known example of singular K3 surfaces \cite{Schutt}. Its transcendental cycles generate a two dimensional Galois representation $V({\mathscr{F}_0})$ that is modular, associated to which is a weight-3 newform of level 16
\begin{equation}
\eta(4z)^6 \in S_3(\Gamma_0(16),\chi_{16}).
\end{equation}
Here the Dirichlet character $\chi_{16}$ is defined by
\begin{equation}
\chi_{16}: (\mathbb{Z}/16\mathbb{Z})^\times \rightarrow \mathbb{C},~\text{with}~\chi_{16}(5)=1,~\chi_{16}(15)=-1.
\end{equation}
Modularity of  $V({\mathscr{F}_0})$ means that its zeta function at a good prime $p$ is of the form
\begin{equation} \label{eq:zetafnFermatQ}
1-b_p(\mathscr{F}_0)T+p^2T^2,
\end{equation}
where $b_p(\mathscr{F}_0)$ is the $p$-th coefficient of the $\textbf{q}$-expansion of the weight-3 newform $\eta(4z)^6$. 

The elliptic partner of the Fermat quartic is the elliptic curve
\begin{equation}
\mathscr{E}_2:~y^2=x(x-1)(x-2),
\end{equation}
whose Weierstrass minimal model is \ref{eq:lmfdb32a3}. Notice that the $j$-invariant of $\mathscr{E}_2$ is 1728, and it admits CM. The zeta function of $H^1_{\text{\'et}}(\mathscr{E}_2,\mathbb{Q}_\ell)$ at a good prime $p$ is of the form
\begin{equation} \label{eq:32aaLMFDB}
1-a_p(\mathscr{E}_2)T+pT^2=(1-\pi^1_p(\mathscr{E}_2)T)(1-\pi^2_p(\mathscr{E}_2)T),
\end{equation}
where $a_p(\mathscr{E}_2)$ is the $p$-th coefficient of the $\textbf{q}$-expansion of the weight-2 newform labeled as \textbf{32.2.a.a} in LMFDB \cite{LMFDB}. The symmetric square of \ref{eq:32aaLMFDB} is of the form 
\begin{equation} \label{eq:SQlambda2}
(1-pT)(1-(a^2_p(\mathscr{E}_2)-2p)T+p^2T^2).
\end{equation}
So one might be wondering what is the relation between the quadratic factor of \ref{eq:SQlambda2} and the zeta function \ref{eq:zetafnFermatQ} of Fermat quartic? In fact, it is very interesting that we have
\begin{equation}
b_p(\mathscr{F}_0)=a^2_p(\mathscr{E}_2)-2p,
\end{equation}
or equivalently 
\begin{equation} \label{eq:sqFermatElliptic}
1-(a^2_p(\mathscr{E}_2)-2p)T+p^2T^2=1-b_p(\mathscr{F}_0)T+p^2T^2.
\end{equation}
Hence we can say the modular form $\eta(4z)^6$ associated to the Fermat quartic arises from the symmetric square of \textbf{32.2.a.a} \cite{Ono}.

\subsection{Deligne's periods for Fermat quartic}

The transcendental cycles of the Fermat quartic defines a two dimensional pure motive $\mathbf{M}_0$ over $\mathbb{Q}$, whose \'etale realization is the two dimensional Galois representation $V(\mathscr{F}_0)$ in Section \ref{sec:ModularFermat}. Thus the $L$-function of $\mathbf{M}_0$ is just the $L$-function associated to the weight-3 newform $\eta(4z)^6$ \cite{Schutt}
\begin{equation}
L(\mathbf{M}_0, s)=L(\eta(4z)^6,s).
\end{equation}
From Mellin transform, $L(\eta(4z)^6,s)$ has an integral representation \cite{Diamond}
\begin{equation} \label{eq:mellintransformL}
L(\eta(4z)^6,s)=\frac{(2 \pi)^s}{\Gamma(s)} \int_0^\infty \eta(4iz)^6z^s \frac{dz}{z}.
\end{equation}
The Hodge realization of $\mathbf{M}_0$ is a two dimensional pure Hodge structure whose Hodge decomposition only has (2,0) and (0,2) parts. Moreover, this pure Hodge structure is completely determined by the expansion \ref{eq:fermatTwoformExpansion} of the holomorphic twoform $\Omega_F$ on the Fermat quartic $\mathscr{F}_0$.

The computation of the Deligne's period $c^+(\mathbf{M}_0)$ for $\mathbf{M}_0$ immediately follows from the method in the paper \cite{YangDeligne}. More explicitly, $c^+(\mathbf{M}_0)$ is given by the pairing of a cohomology cycle of $H^2(X,\mathbb{Q})$ and $\Omega_F$.  From Sections \ref{sec:mirrorfamilyconstruction} and \ref{sec:mirrormapLegendre}, the quotient $\varpi_1(2)/\varpi_0(2)$ is given by
\begin{equation}
\frac{\varpi_1(2)}{\varpi_0(2)}=\frac{-1+i}{2}.
\end{equation}
Hence from the method in \cite{YangDeligne}, we deduce that there exist rational numbers $r_i \in \mathbb{Q}$ such that $c^+(\mathbf{M}_0)$ is of the form
\begin{equation}
c^+(\mathbf{M}_0)=(1+i)\left[ r_0+r_1\frac{-1+i}{2}+r_2 \left(\frac{-1+i}{2} \right)^2 \right] \varpi^2_0(2).
\end{equation}
Since Deligne's period is only well-defined up to a nonzero rational multiple, we immediately learn that there exist two rational numbers $s_1$ and $s_2$ such that 
\begin{equation} \label{eq:Deligneperiods1}
c^+(\mathbf{M}_0)=(s_1+s_2 i) \varpi^2_0(2);~s_i \in \mathbb{Q}.
\end{equation}
But from the construction of Deligne's period, $c^+(\mathbf{M}_0)$ must be a real number \cite{DeligneL}, which uniquely determines the values of $s_1$ and $s_2$ up to a nonzero rational multiple. Similarly, there exist two rational numbers $s_3$ and $s_4$ such that the Deligne's period $c^-(\mathbf{M}_0)$ is given by
\begin{equation} \label{eq:Deligneperiods2}
c^-(\mathbf{M}_0)=(s_3+s_4i) \varpi^2_0(2);~s_i \in \mathbb{Q}.
\end{equation}
From its construction, $c^-(\mathbf{M}_0)$ must be a purely imaginary number \cite{DeligneL}, which uniquely determines the values of $s_3$ and $s_4$ up to a nonzero rational multiple. Furthermore, from formulas \ref{eq:singularpointscoords} and \ref{eq:varpi0theta3}, $\varpi_0(2)$ is equal to the value of $\theta^2_3(0,q)$ at $q=\exp(\pi i(-1+i)/2)$, i.e.
\begin{equation}
\varpi_0(2)=\theta^2_3(0,-i e^{-\pi/2}).
\end{equation}

\subsection{The verification of Deligne's conjecture}

From \cite{DeligneL}, the Tate twist $\mathbf{M}_0 \otimes \mathbb{Q}(n)$ is critical if and only if $n=1,2$. Deligne's conjecture predicts that $c^+(\mathbf{M}_0 \otimes \mathbb{Q}(1))$ (resp. $c^+(\mathbf{M}_0 \otimes \mathbb{Q}(2))$) is a rational multiple of $L(\mathbf{M}_0 \otimes \mathbb{Q}(1),0)$  (resp. $L(\mathbf{M}_0 \otimes \mathbb{Q}(2),0)$) \cite{DeligneL,YangDeligne,YangAttractor}. From \cite{DeligneL,YangDeligne}, we learn that
\begin{equation}
\begin{aligned}
c^+(\mathbf{M}_0 \otimes \mathbb{Q}(1))&=(2 \pi i)c^-(\mathbf{M}_0),\\
c^+(\mathbf{M}_0 \otimes \mathbb{Q}(2))&=(2 \pi i)^2c^+(\mathbf{M}_0).
\end{aligned}
\end{equation}
On the other hand, the $L$-function of a Tate twist is given by \cite{DeligneL,YangDeligne}
\begin{equation}
L(\mathbf{M}_0 \otimes \mathbb{Q}(n),s)=L(\mathbf{M}_0,n+s),
\end{equation}
hence formula \ref{eq:mellintransformL} implies
\begin{equation} \label{eq:LfunctionIntegrals}
\begin{aligned}
L(\mathbf{M}_0 \otimes \mathbb{Q}(1),0)&=2 \pi \int_0^\infty \eta(4iz)^6 dz,\\
L(\mathbf{M}_0 \otimes \mathbb{Q}(2),0)&=(2 \pi)^2\int_0^\infty \eta(4iz)^6 z dz.
\end{aligned}
\end{equation}

Now we will numerically verify that the critical motives $\mathbf{M}_0 \otimes \mathbb{Q}(n)$ with $n=1,2$ satisfy the predictions of Deligne's conjecture. First, the numerical value of $\theta^4_3(0,-i e^{-\pi/2})$ can be evaluated to a very high precision by Mathematica
\begin{equation}
\theta^4_3(0,-i e^{-\pi/2})=- i \,1.3932039296856768591842462603253682426574812175156\cdots,
\end{equation}
which is purely imaginary. Hence in the formulas \ref{eq:Deligneperiods1} and \ref{eq:Deligneperiods2}, we can choose
\begin{equation}
s_1=0,s_2=1,s_3=1,s_4=0,
\end{equation}
i.e. we have 
\begin{equation}
\begin{aligned}
c^+(\mathbf{M}_0 \otimes \mathbb{Q}(1))&=(2 \pi i)\,\theta^4_3(0,-i e^{-\pi/2}),\\
c^+(\mathbf{M}_0 \otimes \mathbb{Q}(2))&=i(2 \pi i)^2\,\theta^4_3(0,-i e^{-\pi/2}).
\end{aligned}
\end{equation}
The integrals in formula \ref{eq:LfunctionIntegrals} can also be numerically evaluated. In this paper, we have computed the first 300 digits of them and here we give the first 50 digits
\begin{equation} 
\begin{aligned}
L(\mathbf{M}_0 \otimes \mathbb{Q}(1),0)&=0.5471099038066191597091924851761161358148431807064 \cdots,\\
L(\mathbf{M}_0 \otimes \mathbb{Q}(2),0)&=0.8593982272525466034362619724763196497376070564774 \cdots.
\end{aligned}
\end{equation}
From these numerical results, we immediately obtain
\begin{equation}
\begin{aligned}
c^+(\mathbf{M}_0 \otimes \mathbb{Q}(1))&=16 \,L(\mathbf{M}_0 \otimes \mathbb{Q}(1),0),\\
c^+(\mathbf{M}_0 \otimes \mathbb{Q}(2))&=-64\, L(\mathbf{M}_0 \otimes \mathbb{Q}(2),0),
\end{aligned}
\end{equation}
which indeed satisfy the predictions of Deligne's conjecture \cite{DeligneL,YangDeligne}.

\section{Conclusions and further prospects} \label{sec:conclusions}

In this paper, we have studied the highly interesting connections between the mirror symmetry of K3 surfaces and the geometry of the Legendre family of elliptic curves. Using the quadratic transformations of hypergeometric functions, we have found interesting relations between the periods of the holomorphic twoform of the Dwork family (Fermat pencil) of K3 surfaces and the periods of the holomorphic oneform of the Legendre family. Then we have shown that the mirror map of the Dwork family is the same as the period map of the Legendre family, which is a crucial result of this paper that provides important insights into the nature of the mirror symmetry of K3 surfaces. For example, it gives an interesting interpretation to the modularity of counting functions for K3 surfaces from the mirror symmetry point of view. Furthermore, these results imply the existence of interesting connections between the arithmetic geometry of the Dwork family and the geometry of the Legendre family, e.g. the Shioda-Inose structures. 

We have also explored the potential relations between the zeta functions of smooth fibers of the Fermat pencil and that of the smooth fibers of the Legendre family. In particular, we have studied the relations between the weight-3 newform $\eta(4z)^6$ associated to the Fermat quartic and the weight-2 newform \textbf{32.2.a.a} associated to the smooth fiber at $\lambda=2$ of the Legendre family. More concretely, $\eta(4z)^6$ can be considered as the symmetric square of \textbf{32.2.a.a}. We have also computed the Deligne's periods of the Fermat quartic, which are given by special values of the theta function $\theta_3$; then numerically we have shown that they satisfy the predictions of Deligne's conjecture.

There are still many open questions left unaddressed, and here we list several interesting ones that come to our mind:
\begin{enumerate}
\item Are there any connections between the results of this paper and the homological mirror symmetry for the quartic K3 surfaces studied in the paper \cite{Seidel}?

\item Could the results in Section \ref{sec:mirrormapLegendre} be applied to study the modularities of counting functions for K3 surfaces?

\item Could the results in this paper provide interesting links between the mirror symmetry of K3 surfaces and that of elliptic curves studied in the paper \cite{DF}? 

\item What is the relation between the arithmetic geometry of the singular fiber $\mathscr{F}_1$ of the Fermat pencil \ref{eq:fermatpencilK3} and the elliptic curve of the Legendre family over the point $\lambda=2 \sqrt{2}-2$, whose $j$-invariant is 8000?

\item Whether the zeta functions of smooth fibers of the Fermat pencil (or Dwork family) can be computed using the zeta functions of smooth fibers of the Legendre family?
\end{enumerate}

\section*{Acknowledgments}

The author is grateful to Brandon Rayhaun for helpful correspondences.

\appendix

\section{A review of Weil conjectures} \label{sec:Weilconjectures}

The concept of zeta functions of a non-singular variety comes from points-counting modulo a prime number. Suppose $X$ is an $n$-dimensional non-singular variety defined over $\mathbb{Q}$, which has an integral model $\mathcal{X}$ defined over $\mathbb{Z}$. Modulo a prime number $p$, $\mathcal{X}$ defines a variety over the finite field $\mathbb{F}_p:=\mathbb{Z}/p\mathbb{Z}$, which will be denoted by $X/\mathbb{F}_p$. We say $p$ is a good prime of $X$ if $X/\mathbb{F}_p$ is non-singular. 

Suppose $p$ is a good prime of $X$ and $m$ is a positive integer. Recall that $\mathbb{F}_{p^m}$ is the unique degree-$m$ extension of $\mathbb{F}_p$. Since $\mathbb{F}_p$ is a subfield of $\mathbb{F}_{p^m}$, the variety $X/\mathbb{F}_p$ is naturally a variety over $\mathbb{F}_{p^m}$. Let $N_m$ be the number of points of $X/\mathbb{F}_p$ with coordinates lie in $\mathbb{F}_{p^m}$. The zeta function $\zeta(X,p,T)$ is by definition the generating series
\begin{equation}
\zeta(X,p,T):=\exp \left( \sum_{m=1}^{\infty} \frac{N_m}{m}  \,T^m     \right)
\end{equation}
A priory, $\zeta(X,p,T)$ is only a formal power series in $T$, but Weil's conjectures claim that $\zeta(X,p,T)$ is in fact a rational function in $T$ that can be expressed as
\begin{equation}
\zeta(X,p,T)=\frac{P_1(X,p,T) \cdots P_{2n-1}(X,p,T)}{P_0(X,p,T) \cdots P_{2n}(X,p,T)},
\end{equation}
where each $P_i(X,p,T)$ is an integral polynomial. Furthermore, $P_0(X,p,T)$ and $P_{2n}(X,p,T)$ are of very simple forms
\begin{equation}
P_0(X,p,T)=1-T,~P_{2n}(X,p,T)=1-p^nT.
\end{equation}
The variety $X$ defines an $n$-dimensional complex manifold $X(\mathbb{C})$, and Weil conjectures claim that
\begin{equation}
\text{deg}\,P_i(X,p,T)=\text{dim}_{\mathbb{Q}}H^i(X(\mathbb{C}),\mathbb{Q}).
\end{equation}
The rationality part of Weil conjectures is first proved by Dwork using $p$-adic analysis \cite{Dwork}.  It can also be proved by the existence of a suitable Weil cohomology theory, e.g. \'etale cohomology theory, and the polynomial $P_i(X,p,T)$ is given by the characteristic polynomial of the (geometric) Frobenius action on the \'etale cohomology group $H^i_{\text{\'et}}(X,\mathbb{Q}_\ell)$ \cite{MilneEC}
\begin{equation}
P_i(X,p,T)=\det \left(\text{Id}-T\,\text{Fr}|_{H^i_{\text{\'et}}(X,\mathbb{Q}_\ell)} \right).
\end{equation}
Over the complex field $\mathbb{C}$, the polynomial $P_i(X,p,T)$ factors into the products of linear polynomials
\begin{equation}
P_i(T)=\prod_j (1- \alpha_{ij}T).
\end{equation}
The `Riemann hypothesis' part of Weil conjectures claims that the absolute value of the algebraic number $\alpha_{ij}$ satisfies
\begin{equation}
|\alpha_{ij}|=p^{i/2},
\end{equation}
which is first proved by Deligne. 

Let us now look at the zeta functions of K3 surfaces. Suppose $X$ is an algebraic K3 surface defined over $\mathbb{Q}$. The \'etale cohomology group $H^2_{\text{\'et}}(X,\mathbb{Q}_\ell)$ is a 22-dimensional representation of the absolute Galois group $\text{Gal}(\overline{\mathbb{Q}}/\mathbb{Q})$. Suppose $p$ is a good prime of $X$, then the zeta function of $X$ at $p$ is of the form
\begin{equation} \label{eq:zetafnsofK3surfaces}
\zeta(X,p,T)=\frac{1}{(1-T)\,P_2(X,p,T)\,(1-p^2\,T)},
\end{equation}
where $P_2(X,p,T)$ is an integral polynomial of degree 22 given by
\begin{equation}
P_2(X,p,T)=\text{det} \left(\text{Id}-T\,\text{Fr}_p|_{H^2_{\text{\'et}}(X,\mathbb{Q}_\ell)} \right).
\end{equation}
The polynomial $P_2(X,p,T)$ can further factorize into the products of lower degrees polynomials. More concretely, $H^2_{\text{\'et}}(X,\mathbb{Q}_\ell)$ splits into the direct sum of two sub-representations
\begin{equation}
H^2_{\text{\'et}}(X,\mathbb{Q}_\ell)=V^a_\ell \oplus V^t_\ell,
\end{equation}
where $V^a_\ell $ is generated by the algebraic cycles of $X$ and $V^t_\ell$ is generated by the transcendental cycles of $X$. Hence $V^a_\ell$ is a $\rho(X)$ dimensional representation of $\text{Gal}(\overline{\mathbb{Q}}/\mathbb{Q})$, while $V^t_\ell$ is a $22-\rho(X)$ dimensional representation of $\text{Gal}(\overline{\mathbb{Q}}/\mathbb{Q})$. The polynomial $P_2(X,p,T)$ factorize into to the product
\begin{equation}
P_2(X,p,T)=P^a_2(X,p,T) P^t_2(X,p,T).
\end{equation}
Here $P^a_2(X,p,T)$ is an integral polynomial with degree $\rho(X)$ given by $V^a_\ell$
\begin{equation}
P^a_2(X,p,T)=\text{det} \left(\text{Id}-T\,\text{Fr}_p|_{V^a_\ell} \right),
\end{equation}
and $P^t_2(X,p,T)$ is an integral polynomial with degree $22-\rho(X)$ given by $V^t_\ell$
\begin{equation}
P^t_2(X,p,T)=\text{det} \left(\text{Id}-T\,\text{Fr}_p|_{V^t_\ell} \right).
\end{equation}

\end{document}